\documentclass{article}%
\usepackage{amsmath}
\usepackage{amsfonts}
\usepackage{amssymb}
\usepackage{graphicx}
\usepackage[singlespacing]{setspace}%
\providecommand{\U}[1]{\protect\rule{.1in}{.1in}}

\begin{document}

\title{Robust and efficient estimation of multivariate scatter and location }
\author{Ricardo A. Maronna$^{1}$ and V\'{\i}ctor J. Yohai$^{2}$\\$^{1}$National University of La Plata, La Plata, Argentina (rmaronna@retina.ar)\\$^{2}$University of Buenos Aires and CONICET, Argentina (victoryohai@gmail.com)}
\date{}
\maketitle

\begin{abstract}
We deal with the equivariant estimation of scatter and location for
$p$-dimensional data, giving emphasis to scatter. It is important that the
estimators possess both a high efficiency for normal data and a high
resistance to outliers, that is, a low bias under contamination. The most
frequently employed estimators are not quite satisfactory in this respect. The
Minimum Volume Ellipsoid (MVE) and Minimum Covariance Determinant (MCD)
estimators are known to have a very low efficiency. S-Estimators with a
monotonic weight function like the bisquare have a low efficiency for small
$p,$ and their efficiency tends to one with increasing $p$. Unfortunately,
this advantage is paid for by a serious loss of robustness for large $p$.

We consider four families of estimators with controllable efficiencies whose
performance for moderate to large $p$ has not been explored to date:
S-estimators with a non-monotonic weight function (Rocke 1996), MM-estimators,
$\tau$-estimators, and the Stahel-Donoho estimator. Two types of starting
estimators are employed: the MVE computed through subsampling, and a
semi-deterministic procedure proposed by Pe\~{n}a and Prieto (2007) for
outlier detection.

A simulation study shows that the Rocke estimator starting from the
Pe\~{n}a-Prieto estimator and with an adequate tuning, can simultaneously
attain high efficiency and high robustness for $p\geq15,$ and the MM estimator
can be recommended for $p$%
$<$%
15.


\end{abstract}

\section{Introduction}

Consider a sample $\mathbf{X=}\{\mathbf{x}_{1},...,\mathbf{x}_{n}\}\subset
R^{p}.$ We look for substitutes $\widehat{\mathbf{\mu}}\in R^{p}$\textbf{
}and\textbf{ }$\widehat{\mathbf{\Sigma}}\in R^{p\times p}$ of the sample mean
vector and covariance matrix, that are resistant to atypical observations. We
also want estimators that have a high efficiency for normal samples$.$ As a
measure of robustness we consider not only the breakdown point but also the
maximum expected Kullback-Leibler divergence between the estimator and the
true value. under contamination.

The most frequently employed estimators are not quite satisfactory in this
respect. The Minimum Volume Ellipsoid (MVE) and the Minimum Covariance
Determinant (MCD) estimators (Rousseeuw 1985) are known to have a low
efficiency. This efficiency can be increased by means of a \textquotedblleft
one-step reweighting\textquotedblright. Croux and Haesbroeck (1999, Tables VII
and VIII) computed the finite-sample efficiencies of the reweighted MCD;
although they are much higher than for the \textquotedblleft
raw\textquotedblright\ estimator, they are still low if one wants a high
breakdown point.  S-Estimators (Davies 1987) with a monotonic weight function
like the bisquare have a low efficiency for small $p.$ Rocke (1996) showed
that their efficiency tends to one with increasing $p$; unfortunately, this
advantage is paid for with a serious loss of robustness for large $p$.

We restrict ourselves to equivariant estimators. There exist many
non-equivariant proposals; but the comparison between equivariant and
non-equivariant estimators is difficult. In particular, a non-equivariant
estimator is more difficult to tune for a given efficiency, since the latter
depends on the correlations.

Among the published equivariant proposals, there are four families of
estimators with controllable efficiencies: non-monotonic S-estimators (Rocke
1996), MM-estimators (Tatsuoka and Tyler 2000), $\tau$-estimators (Lopuhaa
1991) and the estimator proposed independently by Stahel (1981) and Donoho
(1982) but their behavior for large dimensions has not been explored to date.
We compare their behaviors employing different weight functions. A simulation
study shows that the Rocke and MM estimators, with an adequate weight function
and an adequate tuning, can simultaneously attain high efficiency and high robustness.

It will be seen below that if we have a good $\widehat{\Sigma},$ it is easy to
find a good equivariant $\widehat{\mathbf{\mu}},$ but the converse is not
true. For this reason we shall put more emphasis on the estimation of the
scatter matrix.

Since all the considered estimators are based on the iterative minimization of
a non-convex function, the starting values are crucial. Subsampling is the
standard way to compute starting values; but we shall see that a
semi-deterministic equivariant procedure proposed by Pe\~{n}a and Prieto
(2007) may yield both shorter computing times and better statistical performances.

In Section \ref{secMesti} we describe monotonic M-estimators; Section
\ref{secMinScale} deals with estimators based on the minimization of a robust
scale of Mahalanobis distances. Sections \ref{secMM} and \ref{Sec_StaDono}
deal with MM and Stahel-Donoho estimators respectively. In Section
\ref{secRho} we discuss the choice of the $\rho-$function for MM- and $\tau
-$estimators. Section \ref{SecComputing} deals with computational details. In
Section \ref{secSimula} the estimators are compared through a simulation
study. In Section \ref{secReal} the estimators are applied to a real data set.
Finally Section \ref{secConclu} summarizes the results. Section
\ref{SecAppend} is an Appendix containing the full results of the simulations,
the approximations for the tuning constants and some details on the Rocke and
the Pe\~{n}a-Prieto procedures.

\section{Monotonic M-estimators\label{secMesti}}

For $\mathbf{x,\mu}\in R^{p}$ and $\mathbf{\Sigma\in}R^{p\times p}$ define the
(squared) Mahalanobis distance as
\[
d\left(  \mathbf{x,\mu,\Sigma}\right)  =\left(  \mathbf{x-\mu}\right)
^{\prime}\Sigma^{-1}\left(  \mathbf{x-\mu}\right)  .
\]

Let $W$ be a bounded nonincreasing \textquotedblleft weight
function\textquotedblright. Then monotonic M-estimators (Maronna 1976) are
defined as solutions of
\begin{align}
\frac{1}{n}\sum_{i=1}^{n}W\left(  d_{i}\right)  \left(  \mathbf{x-\mu}\right)
\left(  \mathbf{x-\mu}\right)  ^{\prime}  &  =\mathbf{\Sigma~\ \ }%
\label{defMesti}\\
\frac{1}{n}\sum_{i=1}^{n}W\left(  d_{i}\right)  \left(  \mathbf{x-\mu}\right)
&  =0 \label{defMesti2}%
\end{align}
where for brevity we put
\[
\mathrm{\ }d_{i}=d\left(  \mathbf{x}_{i}\mathbf{,\mu,\Sigma}\right)  .
\]
\qquad

The uniqueness of the solutions requires that $W\left(  d\right)  d$ be
nondecreasing. Unfortunately, this implies (Maronna, 1976) that the breakdown
point is $\leq1/\left(  p+1\right)  ,$ which makes these estimators unreliable
except for small $p.$ Besides, this fact holds even if $\mathbf{\mu}$ is
known, while the asymptotic breakdown point of $\widehat{\mathbf{\mu}}$ with
known $\mathbf{\Sigma}$ is 0.5 with an adequate $W.$ This shows that the main
problem to attain high robustness is the scatter matrix.

\section{Estimators based on the minimization of a robust
scale\label{secMinScale}}

For $\mathbf{d=}\left(  d_{1},...,d_{n}\right)  $ let $S\left(  \mathbf{d}%
\right)  $ be a robust scale. Put
\[
\mathbf{d}\left(  \mathbf{\mu,\Sigma}\right)  =\left(  d\left(  \mathbf{x}%
_{1},\mathbf{\mu,\Sigma}\right)  ,...,d\left(  \mathbf{x}_{1},\mathbf{\mu
,\Sigma}\right)  \right)  .
\]
A general family of estimators can be defined by
\begin{equation}
\left(  \widehat{\mathbf{\mu}}\mathbf{,}\widetilde{\mathbf{\Sigma}}\right)
=\arg\min S\left(  \mathbf{d}\left(  \mathbf{\mu,\Sigma}\right)  \right)
,\ \mathbf{\mu\in}R^{p},\ \mathbf{\Sigma\in}R^{p},\ |\mathbf{\Sigma}|=1,
\label{defEscaMini}%
\end{equation}
where the condition $|\mathbf{\Sigma}|=1$ rules out trivial solutions with
$\mathbf{\Sigma\rightarrow\infty}$.

If $S\left(  \mathbf{d}\right)  =\mathrm{Median}\left(  \mathbf{d}\right)  $
we have the \textquotedblleft Minimum Volume Ellipsoid\textquotedblright%
\ (MVE) estimator, and if $S$ is a trimmed mean, we have the \textquotedblleft
Minimum Covariance Determinant\textquotedblright\ (MCD) estimator, both
proposed by Rousseeuw (1985). The first one is very robust, but has a null
asymptotic efficiency; the second is very popular, but its asymptotic
efficiency is very low; see (Paindaveine and Van Bever, 2014) and references
therein, and its maximum contamination bias increases rapidly with $p$
(Agostinelli et al, 2015, Table 1).

The condition $|\mathbf{\Sigma}|=1$ means that we estimate the
\textquotedblleft shape\textquotedblright\ of the scatter. Given the shape,
the \textquotedblright size\textquotedblright\ can easily estimated to yield
consistency at the normal model (Maronna et al., Section 6.3.2). A simple way
is to put%
\begin{equation}
\widehat{\mathbf{\Sigma}}\mathbf{=}\frac{\mathrm{Median}\left(  \mathbf{d}%
\left(  \widehat{\mathbf{\mu}}\mathbf{,}\widetilde{\mathbf{\Sigma}}\right)
\right)  }{\mathrm{Median}\left(  \chi_{p}^{2}\right)  }\widetilde
{\mathbf{\Sigma}}. \label{CorreTama}%
\end{equation}
Instead of the median, one could use more efficient scales, such as an
M-scale, but exploratory simulations indicate that they do not yield better results.

\subsection{S-estimators\label{secSestim}}

Let $S=S\left(  d_{1},..,d_{n}\right)  $ be a scale M-estimator defined as
solution of%

\begin{equation}
\frac{1}{n}\sum_{i=1}^{n}\rho\left(  \frac{d_{i}}{S}\right)  =\delta,
\label{defMScale}%
\end{equation}
where $\delta\in\left(  0,1\right)  $ controls the breakdown point, and
$\rho\left(  t\right)  \in\lbrack0,1]$ is smooth and nondecreasing in
$t\geq0,$ with $\rho\left(  0\right)  =0$ and $\max\rho=1.$ Then S-estimators
(Davies 1987) are defined as solutions of (\ref{defEscaMini}) with $S$ given
by (\ref{defMScale}).

The maximum finite-sample replacement breakdown point is attained when%

\begin{equation}
\delta=0.5\left(  1-\frac{p}{n}\right)  , \label{BDP-S}%
\end{equation}
and its value is equal to this $\delta.$ See (Maronna et al., 2006, Section 6.4.2).

A popular $\rho$ is the bisquare given by
\begin{equation}
\rho\left(  d\right)  =\left\{
\begin{array}
[c]{ccc}%
1-\left(  1-d\right)  ^{3} & \mathrm{if} & d\leq1\\
1 & \mathrm{if} & d>1.
\end{array}
\right.  \label{defiBis}%
\end{equation}
Note that the usual bisquare $\rho$ employed for regression is actually
$\rho_{\mathrm{bis}}\left(  t\right)  =\rho\left(  t^{2}\right)  .$ However,
since we are dealing with the \emph{squared} distances, we employ in
(\ref{defiBis}) $\rho_{\mathrm{bis}}\left(  \sqrt{d}\right)  =\rho\left(
d\right)  .$

It is easy to show that S-estimators satisfy the \textquotedblleft estimating
equations\textquotedblright%
\begin{align}
\frac{1}{n}\sum_{i=1}^{n}W\left(  \frac{d_{i}}{S}\right)  \left(
\mathbf{x-\mu}\right)  \left(  \mathbf{x-\mu}\right)  ^{\prime}  &
=\mathbf{\Sigma~\ \ }\label{EstimEquat_SE}\\
\frac{1}{n}\sum_{i=1}^{n}W\left(  \frac{d_{i}}{S}\right)  \left(
\mathbf{x-\mu}\right)   &  =0\label{EstimEquat_SE2}\\
\frac{1}{n}\sum_{i=1}^{n}\rho\left(  \frac{d_{i}}{S}\right)   &
=\delta\label{EstimEquat_SE3}%
\end{align}
with $W=\rho^{\prime}.$ That is, they satisfy the equations (\ref{defMesti}%
)-(\ref{defMesti2}) which define monotonic M-estimators, with weight function
$W=\rho^{\prime}.$ Here, since $\rho$ is bounded $W\left(  d\right)  d$ is not
a nondecreasing function, and therefore this case is different from
monotonic\ M-estimators. In particular, the breakdown point is not bounded by
$\left(  1+p\right)  ^{-1};$ as shown by (\ref{BDP-S}).

For the bisquare, the weight function is%
\[
W(t)=3\left(  1-t\right)  ^{2}\mathrm{I}\left(  t\leq1\right)
\]
(where $\mathrm{I}\left(  .\right)  $ denotes the indicator), which is
decreasing. It seems intuitive that the weights of the observations should
decrease with their \textquotedblleft outlyingness\textquotedblright. However
it will be seen in the next Section that monotonicity is not necessarily favorable.

\subsection{S-estimators with a non-monotonic weight
function\label{secNonMono}}

Rocke (1996) showed that if $W$ is nonincreasing, the efficiency of the
estimator tends to one when $p\rightarrow\infty.$ A similar result was derived
by Kent and Tyler (1996, page 1363) for their constrained M-estimators.

Table \ref{TabEfiSbis} shows the efficiencies (to be defined later) of the
bisquare S-estimator of scatter for normal $p$-dimensional data.

\begin{center}%
\begin{table}[tbp] \centering
\begin{tabular}
[c]{cccccccc}\hline
$p$ & 2 & 5 & 10 & 20 & 30 & 40 & 50\\
Efficiency & 0.427 & 0.793 & 0.930 & 0.976 & 0.984 & 0.990 & 0.992\\\hline
\end{tabular}
\caption{Efficiencies of the S-estimator with bisquare weights for dimension  $p$ }\label{TabEfiSbis}%
\end{table}%

\end{center}

However, it will be seen that the price for this increase in efficiency is a
decrease in robustness. More precisely, although the breakdown point does not
tend to zero with increasing $p,$ the bias caused by contamination grows
rapidly with $p.$ This fact suggests that we need estimators with a
controllable efficiency. But while in regression the efficiency has to be
controlled to make it higher, here we need to prevent it from becoming
\textquotedblleft too high\textquotedblright.

Based on the fact that for large $p$ the $p$-variate standard normal
distribution $\mathrm{N}_{p}\left(  \mathbf{0,I}\right)  $ is concentrated
\textquotedblleft near\textquotedblright\ the spherical shell with radius
$\sqrt{p},$ Rocke (1996) proposed estimators with non-monotonic weight
functions. Maronna et al. (2006) proposed a modification of Rocke's
\textquotedblleft biflat\textquotedblright\ function, namely%
\begin{equation}
W\left(  d\right)  =\left[  1-\left(  \frac{d-1}{\gamma}\right)  ^{2}\right]
\mathrm{I}\left(  1-\gamma\leq d\leq1+\gamma\right)  \label{RockeWeight}%
\end{equation}
with%
\begin{equation}
\gamma=\min\left(  1,\frac{\chi_{p}^{2}\left(  1-\alpha\right)  }{p}-1\right)
, \label{gama-alfa}%
\end{equation}
where $\chi_{p}^{2}\left(  \beta\right)  $ is the $\beta$-quantile of the
$\chi^{2}$ distribution with $p$ degrees of freedom, and $\alpha$ is
\textquotedblleft small\textquotedblright\ to control the efficiency.

Maronna et al (2006, Sec. 6.8) dealt only with location. The performance of
the respective scatter matrix will be studied below.

\subsection{$\tau-$estimators\label{secTau}}

$\tau-$estimators were proposed by Yohai and Zamar (1988) to obtain robust
regression estimators with controllable efficiency, and later Lopuha\"{a}
(1991) employed the same approach for multivariate estimation. This approach
requires two functions $\rho_{1}$ and $\rho_{2}.$ For given $\left(
\mathbf{\mu,\Sigma}\right)  $ call $\sigma_{0}\left(  \mathbf{\mu,\Sigma
}\right)  $ the solution of%

\[
\frac{1}{n}\sum_{i=1}^{n}\rho_{1}\left(  \frac{d\left(  \mathbf{x}%
_{i},\mathbf{\mu,\Sigma}\right)  }{\sigma_{0}}\right)  =\delta.
\]

\bigskip Then the estimator minimizes the \textquotedblleft$\tau
$-scale\textquotedblright%
\[
\sigma\left(  \mathbf{\mu,\Sigma}\right)  =\sigma_{0}\left(  \mathbf{\mu
,\Sigma}\right)  \frac{1}{n}\sum_{i=1}^{n}\rho_{2}\left(  \frac{d\left(
\mathbf{x}_{i},\mathbf{\mu,\Sigma}\right)  }{\sigma_{0}\left(  \mathbf{\mu
,\Sigma}\right)  }\right)  .
\]

Here
\begin{equation}
\rho_{2}\left(  t\right)  =\rho_{1}\left(  \frac{t}{c}\right)  \label{rho1-2}%
\end{equation}
where $c$ is chosen to regulate the efficiency.

Originally, $\tau$-estimators were proposed to obtain estimators with higher
efficiency than S-estimators for small $p,$ which required $c>1;$ but for
large $p$ we need $c<1$ in order to decrease the efficiency.

\section{MM-estimators\label{secMM}}

MM-estimators were initially proposed by Yohai (1987) to obtain regression
estimators with a controllable efficiency. This approach has been used in the
multivariate setting by Lopuha\"{a} (1992) and Tatsuoka and Tyler (2000). Here
we give a simplified version of the latter.

Let $\left(  \widehat{\mathbf{\mu}}_{0}\mathbf{,}\widehat{\mathbf{\Sigma}}%
_{0}\right)  $ be an initial very robust although possibly inefficient
estimator. Put%
\[
d_{i}^{0}=d\left(  \mathbf{x}_{i}\mathbf{,}\widehat{\mathbf{\mu}}%
_{0}\mathbf{,}\widehat{\mathbf{\Sigma}}_{0}\right)
\]
and call $S$ the respective M-scale%
\begin{equation}
\frac{1}{n}\sum_{i=1}^{n}\rho\left(  \frac{d_{i}^{0}}{S}\right)  =\delta.
\label{S-inicialMM}%
\end{equation}

The estimator is defined by $\left(  \widehat{\mathbf{\mu}},\widehat
{\mathbf{\Sigma}}\right)  $ with $|\widehat{\mathbf{\Sigma}}|=1$ such that%
\begin{equation}
\sum_{i=1}^{n}\rho\left(  \frac{d_{i}}{cS}\right)  =\min, \label{defMM_min}%
\end{equation}
where $d_{i}=d\left(  \mathbf{x}_{i},\widehat{\mathbf{\mu}},\widehat
{\mathbf{\Sigma}}\right)  $ and the constant $c$ is chosen to control efficiency.

It can be shown that the solution satisfies the equations
\begin{align}
\frac{1}{n}\sum_{i=1}^{n}W\left(  \frac{d_{i}}{cS}\right)  \left(
\mathbf{x-\mu}\right)  \left(  \mathbf{x-\mu}\right)  ^{\prime}  &
=\mathbf{\Sigma~\ \ }\label{defMM}\\
\frac{1}{n}\sum_{i=1}^{n}W\left(  \frac{d_{i}}{cS}\right)  \left(
\mathbf{x-\mu}\right)   &  =0\nonumber
\end{align}
with $W=\rho^{\prime},$

Actually, it is not necessary to obtain the absolute minimum in
(\ref{defMM_min}). As with regression MM-estimators (Yohai 1987) it is
possible to show that any solution of (\ref{defMM}) for with the objective
function (\ref{defMM_min}) is lower that for the initial estimator, has the
same asymptotic behavior as the absolute minimum and has a breakdown point at
least as high as the initial estimator.

Like $\tau$-estimators, MM estimators were originally proposed to obtain
estimators with higher efficiency than S-estimators for small $p;$ but here
for large $p$ the constant has to chosen to prevent the efficiency becoming
too high.

\section{The Stahel-Donoho estimator\label{Sec_StaDono}}

Let $M\left(  .\right)  $ and $S\left(  .\right)  $ be univariate location and
dispersion statistics, e.g., the median and MAD. Define for any $\mathbf{y}\in
R^{p}$ the \textit{outlyingness }$\mathit{r}$:
\begin{equation}
r({\mathbf{y}})=\max_{\mathbf{a}}\frac{{|{\mathbf{a}^{\prime}}}\mathbf{y}%
-M({\mathbf{a}^{\prime}\mathbf{X}})|}{S{({\mathbf{a}^{\prime}\mathbf{X}})}},
\label{Outlying}%
\end{equation}
where the supremum is over $\mathbf{a}\in R^{p}$ with $\mathbf{a}%
\neq\mathbf{0}$ or equivalently over the spherical surface $S_{p}%
=\{\mathbf{a}\in R^{p}:\parallel\mathbf{a}\parallel=1\}$. Here $\mathbf{a}%
^{\prime}\mathbf{X}$ denotes $\mathbf{a}^{\prime}\mathbf{x}_{1}%
,\textellipsis,\mathbf{a}^{\prime}\mathbf{x}_{n}$. Let $W$ (the \textit{weight
function}) be a positive function. The Stahel---Donoho estimator of location
and scatter, $(\mathbf{t}(\mathbf{X}),\mathbf{V}(\mathbf{X})),$ is a weighted
mean and covariance matrix, with weights $w_{i}=W(r(\mathbf{x}_{i}))$.

If $W$ is continuous, and $W(r)$ and $W(r)r^{2}$ are bounded for $r\geq0$, the
estimators have asymptotic breakdown point $0.5$ for all $p$ at continuous
multivariate models, if $M$ and $S$ have asymptotic breakdown point $0.5$ (see
Hampel et al.1986). The finite-sample breakdown point was derived by Tyler (1994).

Maronna and Yohai (1995) showed that these estimators have order $\sqrt{n}%
$-consistency. Their asymptotic distribution was given by Zuo et al. (2004).

Maronna and Yohai (1995) recommended a \textquotedblleft
Huber-type\textquotedblright\ $W;$ however, further exploratory simulations
indicate that better results are obtained with the weight function described
in the next section.

The numerical computation of these estimators is difficult. Stahel (1981)
proposed an approximate algorithm based on subsampling, the cost of which
increases rapidly with $p$. Pe\~{n}a and Prieto (2007) proposed a fast
algorithm for outlier detection which combines the projections on a set of
$2p$ deterministic directions that are extrema of the kurtosis, and a set of
random directions. Although this method was originally meant for data
analysis, it offers two further uses. First, the resulting projections can be
employed to compute the Stahel-Donoho estimator; second, the method yields a
robust (but probably inefficient) estimator that can be used as a starting
point for the iterative computing of the estimators described above. Further
details about this procedure are given in Section \ref{SecComputing}

\section{Choosing $\rho$ for MM- and $\tau-$estimators\label{secRho}}

The most popular $\rho$ in robust methods seems to be the bisquare. Yohai and
Zamar (1997) proposed a $\rho$ for regression with certain optimality
properties. A simplified variant of this function is given by Muler et al
(2002). Its version for multivariate estimation has weight function

\begin{equation}
W_{\mathrm{opt}}\left(  d\right)  =\left\{
\begin{array}
[c]{ccc}%
1 & \mathrm{if} & d\leq4\\
q\left(  d\right)  & \mathrm{if} & 4<d\leq9\\
0 & \mathrm{if} & d>9
\end{array}
\right.  , \label{OptiWeight}%
\end{equation}
where%
\[
q\left(  d\right)  =-1.944+1.728d-0.312d^{2}+0.016d^{3}%
\]
is such that $W$ is continuous and differentiable at $d=4$ and $d=9.$ The
respective $\rho$ function is%
\[
\rho\left(  d\right)  =\frac{1}{6.494}\left\{
\begin{array}
[c]{ccc}%
d & \mathrm{if} & d\leq4\\
s\left(  d\right)  & \mathrm{if} & 4<d\leq9\\
6.494 & \mathrm{if} & d>9
\end{array}
\right.  ,
\]
where%
\[
s\left(  d\right)  =3.534-1.944d+0.864d^{2}-0.104d^{3}+0.004d^{4}.
\]

\begin{figure}
[tbh]
\begin{center}
\includegraphics[
height=10.546cm,
width=14.0232cm
]%
{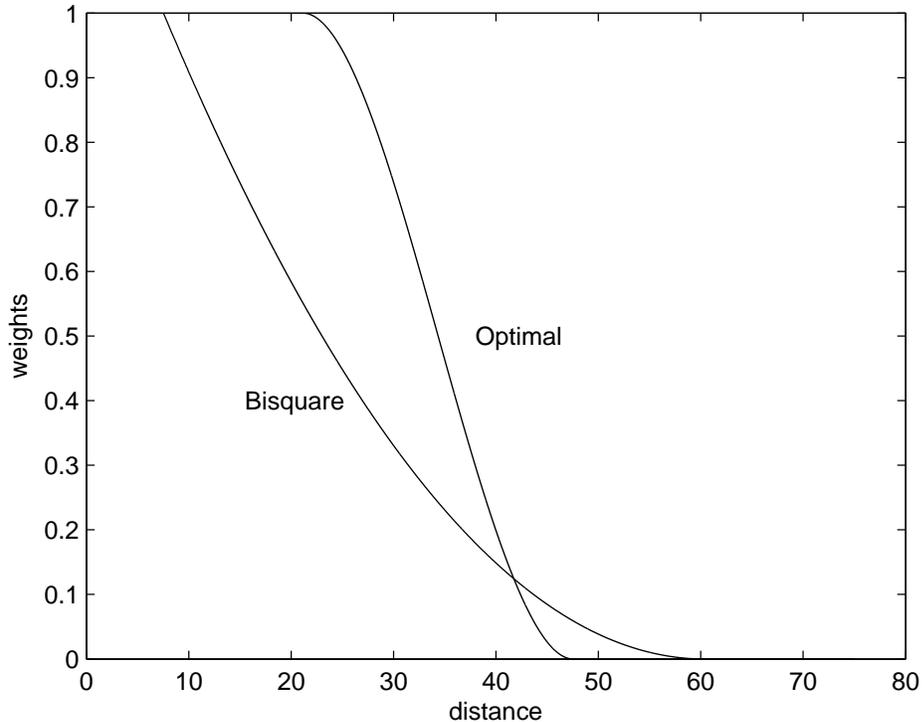}%
\caption{Bisquare and \textquotedblleft optimal\textquotedblright\ weight
functions.}%
\label{figWeigths}%
\end{center}
\end{figure}

Figure \ref{figWeigths} shows the bisquare and \textquotedblleft
optimal\textquotedblright\ weight functions, scaled with their respective
tuning constants for the MM-estimator with 90\% efficiency and $p=30$. It is
seen that the \textquotedblleft optimal\textquotedblright\ $\rho$ yields a
smaller cutoff point.%

\section{Computing issues\label{SecComputing}}

All estimators described above are computed as iterative reweighted means and
covariances, starting from an initial estimator. For S-, $\tau-$ and MM
estimators this algorithm ensures that the objective function descends at each
iteration. This need not happen with the Rocke estimator, which has a
non-monotonic weight function. Maronna et al. (2006, Section 6.4.4) describe
an algorithm which ensures attaining a local minimum.

The (approximate) MVE is computed with 1000 subsamples and using the
improvement described in (Maronna et al., 2006, Section 6.7.3).

Since in all cases we attempt to minimize a non-convex function, the initial
estimator is an essential part of the procedure. The standard way to obtain a
robust and equivariant starting point is subsampling. However, ensuring a high
enough breakdown point with large $p$ may require an impractically large
number of subsamples. Besides, our experiments indicate that the breakdown
point may be much lower than expected when $n/p$ is \textquotedblleft
small\textquotedblright\ (say, $\leq5),$ which is not uncommon with
high-dimensional data sets. For these reasons we need a faster and more
reliable starting point.

Pe\~{n}a and Prieto (2007) proposed an equivariant and semi-deterministic
procedure for outlier detection, based on finding directions that maximize or
minimize the kurtosis of the respective projections, plus a set of random
\textquotedblleft`specific directions\textquotedblright\ aimed at detecting
outliers. Here we employ this procedure (which they call \textquotedblleft
kurtosis plus specific directions\textquotedblright, henceforth abbreviated as
\textquotedblleft KSD\textquotedblright) as an estimator by itself. In the
present setting it would not be competitive with the other estimators because
its efficiency cannot be tuned (see Table \ref{tabEfiKurto} below), but we
shall use it as an initial estimator competing with the sampling-based MVE.

There are no theoretical results on the breakdown point of KSD. However, the
simulations in (Pe\~{n}a and Prieto 2007, table 4) suggest that it can yield
reliable results even with 40\% of outliers. A limited theoretical result is
given in Section \ref{SecAppend}.

\subsection{Computing the Rocke Estimator\label{SecCompuRocke}}

Given a starting point, the Rocke estimator is computed iteratively as
described in Section 9.6.3 of (Maronna et al., 2006).

The form of the weight function ensures that for normal data sets, most of the
data have positive weights. Since real data are seldom normal, it may happen
that for data sets with large $p$ and low ratio $n/p$ the proportion of data
with positive weights is small. If the data set is nearly collinear, this may
cause $\widehat{\mathbf{\Sigma}}$ to be ill-conditioned, which affects the
computation of Mahalanobis distances. For this reason, if at the first
iteration the number of data with positive weights is less than $2p,$ the
tuning constant is enlarged until this number is $\geq2p.$

\section{Simulation\label{secSimula}}

As a reference distribution we take the $p$-variate normal $\mathrm{N}%
_{p}\left(  \mathbf{\mu}_{0},\mathbf{\Sigma}_{0}\right)  .$ In order to
measure the performance of a given estimator $\left(  \widehat{\mathbf{\mu}%
}\mathbf{,}\widehat{\mathbf{\Sigma}}\right)  $ we need a measure of
\textquotedblleft distance\textquotedblright\ between an estimator and the
true value. Recall that the Kullback-Leibler divergence between densities
$f_{1}$ and $f_{2}$ is%
\[
d_{\mathrm{KL}}\left(  f_{1},f_{2}\right)  =\int_{-\infty}^{\infty}\log\left(
\frac{f_{1}\left(  \mathbf{z}\right)  }{f_{2}\left(  \mathbf{z}\right)
}\right)  f_{1}\left(  \mathbf{z}\right)  d\mathbf{z.}%
\]

If both densities belong to the same parametric model with parameter vector
$\mathbf{\theta}$:\ $f_{j}\left(  \mathbf{z}\right)  =f\left(
\mathbf{z,\theta}_{j}\right)  ,$ then $d_{\mathrm{KL}}$ induces a
\textquotedblleft distance\textquotedblright\ between parameters:%

\[
D\left(  \mathbf{\theta}_{1},\mathbf{\theta}_{2}\right)  =d_{\mathrm{KL}%
}\left(  f\left(  .,\mathbf{\theta}_{1}\right)  ,f\left(  .,\mathbf{\theta
}_{1}\right)  \right)  .
\]

In the normal family, for $\mathbf{\mu}$ with known $\Sigma$ we have%
\begin{equation}
D=\left(  \widehat{\mathbf{\mu}}-\mathbf{\mu}_{0}\right)  ^{\prime
}\mathbf{\Sigma}_{0}^{-1}\left(  \widehat{\mathbf{\mu}}-\mathbf{\mu}%
_{0}\right)  , \label{dKL(mu)}%
\end{equation}

and for $\Sigma$ with known $\mu$ we have%
\begin{equation}
D=\mathrm{trace}\left(  \mathbf{\Sigma}_{0}^{-1}\widehat{\mathbf{\Sigma}%
}\right)  -\log|\mathbf{\Sigma}_{0}^{-1}\widehat{\mathbf{\Sigma}}|-p
\label{dKL(Sigma)}%
\end{equation}

Since all estimators are equivariant we may in the simulations take without
loss of generality $\left(  \mathbf{\mu}_{0}\mathbf{,\Sigma}_{0}\right)
=\left(  \mathbf{0,I}\right)  $.

Each estimator is evaluated by $\overline{D}=$ Monte Carlo average of the
Kullback-Leibler divergences $D$ given in (\ref{dKL(mu)})-(\ref{dKL(Sigma)}).

We generate $N=500$ samples $\mathbf{X=[x}_{ij}\mathbf{]}$ of size $n$ from
$N_{p}\left(  \mathbf{0,I}\right)  .$

The estimators compared are:

\begin{itemize}
\item \ Rocke with tuning constant $\alpha;$ see (\ref{gama-alfa})

\item MM with bisquare and \textquotedblleft optimal\textquotedblright%
\ $\rho,$ with tuning constant $c$; see (\ref{defMM})

\item $\tau$ with bisquare and \textquotedblleft optimal\textquotedblright%
\ $\rho,$ with tuning constant $c$; see(\ref{rho1-2})

\item Stahel-Donoho with weight function $W(r)=W_{\mathrm{opt}}\left(
r/c\right)  $ where $W_{\mathrm{opt}}$ is defined in (\ref{OptiWeight})
\end{itemize}

For all estimators we employed both the MVE and KSD estimators as starting
values. The tuning constants were chosen to attain an efficiency of 0.9 (see below).

We add for completeness four other estimators with uncontrollable efficiency:

\begin{itemize}
\item The S-estimator (S-E) with $\delta=0.5$ in (\ref{defMScale}) and
biweight $\rho.$ The \textquotedblleft optimal\textquotedblright\ $\rho$
yielded similar results.

\item The MVE and KSD estimators.

\item The MCD with breakdown point (\ref{BDP-S}) and one-step reweighting,
computed with the code in the library LIBRA (Verboven and Hubert 2005). 
\end{itemize}

All scatter estimators are corrected for \textquotedblleft
size\textquotedblright\ by means of (\ref{CorreTama}), except for the MCD for
which the code applies a consistency correction.

\subsection{No contamination}

Call $\mathbf{C}$ the sample covariance matrix. For each estimator
$\widehat{\mathbf{\Sigma}}$ we define%
\[
\mathrm{efficiency=}\frac{\overline{D}\left(  \mathbf{C}\right)  }%
{\overline{D}\left(  \widehat{\mathbf{\Sigma}}\right)  }.
\]

The constants for each estimator are chosen to attain finite-sample
efficiencies of 0.90. To this end we computed for each estimator its tuning
constants for $n=Kp$ with $K=5,$ 10 and 20 and $p$ between 5 and 50, and then
fitted the constants as functions of $n$ and $p$.

The simulation showed the efficiency cannot be controlled in all cases, namely

\begin{itemize}
\item For $p=15$ the maximum efficiency of the Rocke estimator is 0.876 for
all $\alpha$s, and is still lower for smaller $p$. The explanation is that
when $\alpha$ tends to zero, the estimator does not tend to the covariance
matrix unless $p$ is large enough.

\item The minimum efficiency of the $\tau$-estimators over all constants $c$
tends to one with increasing $p,$ for both $\rho-$functions. In particular, it
is
$>$%
0.95 for $p\geq50.$ The reason is that when $c$ is small, the $\tau-$scale
approaches the M-scale, and therefore the $\tau-$estimators approaches the S-estimator.
\end{itemize}

Table \ref{tabEfiKurto} shows the efficiencies of the KSD estimator.

\begin{center}%
\begin{table}[tbp] \centering
\begin{tabular}
[c]{llll}\hline
$p$ & $n$ & Scatter & Location\\\hline
\multicolumn{1}{r}{10} & \multicolumn{1}{r}{50} & \multicolumn{1}{r}{0.40} &
\multicolumn{1}{r}{0.62}\\
\multicolumn{1}{r}{} & \multicolumn{1}{r}{100} & \multicolumn{1}{r}{0.70} &
\multicolumn{1}{r}{0.85}\\
\multicolumn{1}{r}{} & \multicolumn{1}{r}{200} & \multicolumn{1}{r}{0.86} &
\multicolumn{1}{r}{0.95}\\
\multicolumn{1}{r}{20} & \multicolumn{1}{r}{100} & \multicolumn{1}{r}{0.44} &
\multicolumn{1}{r}{0.62}\\
\multicolumn{1}{r}{} & \multicolumn{1}{r}{200} & \multicolumn{1}{r}{0.80} &
\multicolumn{1}{r}{0.89}\\
\multicolumn{1}{r}{} & \multicolumn{1}{r}{400} & \multicolumn{1}{r}{0.90} &
\multicolumn{1}{r}{0.95}\\
\multicolumn{1}{r}{50} & \multicolumn{1}{r}{250} & \multicolumn{1}{r}{0.47} &
\multicolumn{1}{r}{0.58}\\
\multicolumn{1}{r}{} & \multicolumn{1}{r}{500} & \multicolumn{1}{r}{0.82} &
\multicolumn{1}{r}{0.85}\\
\multicolumn{1}{r}{} & \multicolumn{1}{r}{1000} & \multicolumn{1}{r}{0.93} &
\multicolumn{1}{r}{0.96}\\\hline
\end{tabular}
\caption{Efficiencies of the KSD estimator}\label{tabEfiKurto}%
\end{table}%

\end{center}

It seen that the efficiency depends heavily on the ratio $n/p$ and can be
rather low for $n/p=5.$

\subsection{Contamination}

We deal first with shift contamination. For contamination rate $\varepsilon,$
let $m=[n\varepsilon].$ Given $K,$ we replace the first coordinate:%
\[
x_{i1}\longleftarrow\gamma x_{i1}+K,\ i=1,...,m
\]

The outlier size $K$ is varied between 1 and 12 in order to find the maximum
$\overline{D}.$ The constant $\gamma$ determines the scatter of the outliers.
We employed the values $\varepsilon=0.1$ and 0.2, and $\gamma=0$ and 0.5.

The simulations were run for $p=$ 5,~10, 15, 20 and 30, and $n=mp$ with $m=5,$
10 and 20. Since the complete results are rather bulky, they are given in
Section \ref{secTablas}. Here we give the most important conclusions from
them. Examination of the tables shows that

\begin{itemize}
\item The price paid for the high efficiency of S-E is a large loss of robustness.

\item KSD is always better than MVE as a starting estimator for MM and $\tau$.

\item KSD is generally better than subsampling for S-D.

\item The \textquotedblleft optimal\textquotedblright\ $\rho$ is always better
than the bisquare $\rho$ for both MM and $\tau$

\item In all situations, the best estimators are MM and $\tau$ with
\textquotedblleft optimal\textquotedblright\ $\rho,$ Rocke, and S-D, all
starting from KSD.

\item Although the results for $\gamma=0$ and 0.5 are different, the
comparisons among estimators are almost the same.

\item The relative performances of the estimators for location and scatter are similar.

\item The relative performances of the estimators for $n=5p,$ $10p$ and $20p$
are similar.
\end{itemize}

For these reasons we give in Table \ref{TabResumen} a reduced version of the
results, for $n=10p$ and $\gamma=0,$ and the maximum $\overline{D}$s of the
scatter estimators corresponding to MM and $\tau$ (both with \textquotedblleft
optimal\textquotedblright\ $\rho),$ Rocke and S-D, all starting from KSD. For
completion we add S-E with KSD start, and the reweighted MCD.

The results for estimators with  efficiency less than 0.9 are shown in
italics. 

\begin{center}%
\begin{table}[tbp] \centering
\begin{tabular}
[c]{cccccccc}\hline
$p$ & $\varepsilon$ & MM & $\tau$ & Rocke & S-D & S-E & MCD\\
5 & 0.1 & \multicolumn{1}{r}{0.85} & \multicolumn{1}{r}{0.89} &
\multicolumn{1}{r}{\emph{1.26}} & \multicolumn{1}{r}{0.99} & \emph{1.09} &
\emph{1.99}\\
& 0.2 & \multicolumn{1}{r}{2.27} & \multicolumn{1}{r}{2.46} &
\multicolumn{1}{r}{\emph{3.74}} & \multicolumn{1}{r}{4.53} & \emph{4.38} &
\emph{17.58}\\
\multicolumn{1}{r}{10} & \multicolumn{1}{r}{0.1} & \multicolumn{1}{r}{1.67} &
\multicolumn{1}{r}{1.77} & \multicolumn{1}{r}{\emph{1.53}} &
\multicolumn{1}{r}{1.61} & \emph{3.54} & \emph{6.66}\\
\multicolumn{1}{r}{} & \multicolumn{1}{r}{0.2} & \multicolumn{1}{r}{3.88} &
\multicolumn{1}{r}{4.53} & \multicolumn{1}{r}{\emph{1.67}} &
\multicolumn{1}{r}{7.94} & \emph{11.26} & \emph{21.89}\\
\multicolumn{1}{r}{15} & \multicolumn{1}{r}{0.1} & \multicolumn{1}{r}{2.38} &
\multicolumn{1}{r}{2.98} & \multicolumn{1}{r}{1.95} & \multicolumn{1}{r}{2.26}
& \emph{6.68} & \emph{12.53}\\
\multicolumn{1}{r}{} & \multicolumn{1}{r}{0.2} & \multicolumn{1}{r}{5.68} &
\multicolumn{1}{r}{7.85} & \multicolumn{1}{r}{4.47} &
\multicolumn{1}{r}{12.31} & \emph{19.82} & \emph{28.33}\\
\multicolumn{1}{r}{20} & \multicolumn{1}{r}{0.1} & \multicolumn{1}{r}{3.32} &
\multicolumn{1}{r}{4.59} & \multicolumn{1}{r}{2.49} & \multicolumn{1}{r}{3.00}
& \emph{10.03} & \emph{16.46}\\
\multicolumn{1}{r}{} & \multicolumn{1}{r}{0.2} & \multicolumn{1}{r}{7.90} &
\multicolumn{1}{r}{12.62} & \multicolumn{1}{r}{3.17} &
\multicolumn{1}{r}{17.09} & \emph{25.41} & \emph{32.04}\\
\multicolumn{1}{r}{30} & \multicolumn{1}{r}{0.1} & \multicolumn{1}{r}{5.34} &
\multicolumn{1}{r}{8.56} & \multicolumn{1}{r}{3.03} & \multicolumn{1}{r}{4.64}
& \emph{18.39} & \emph{17.66}\\
\multicolumn{1}{r}{} & \multicolumn{1}{r}{0.2} & \multicolumn{1}{r}{14.21} &
\multicolumn{1}{r}{20.71} & \multicolumn{1}{r}{5.61} &
\multicolumn{1}{r}{29.66} & \emph{49.14} & \emph{34.02}\\\hline
\end{tabular}
\caption{Maximum mean Ds of scatter matrices, for  $n=10p$ and $\gamma=0$. All estimators start from KSD. MM and $\tau$ use ``optimal'' $\rho$.}\label{TabResumen}%
\end{table}%

\end{center}

It is seen that

\begin{itemize}
\item The performance of S-D is competitive for $\varepsilon=0.1,$ but is poor
for $\varepsilon=0.2.$

\item For $p\leq10,$ MM has the best overall performance.

\item For $p\geq15,$ Rocke has the best overall performance.

\item The MCD has a poor performance.
\end{itemize}

\begin{figure}
[tbh]
\begin{center}
\includegraphics[
height=11.2335cm,
width=14.937cm
]%
{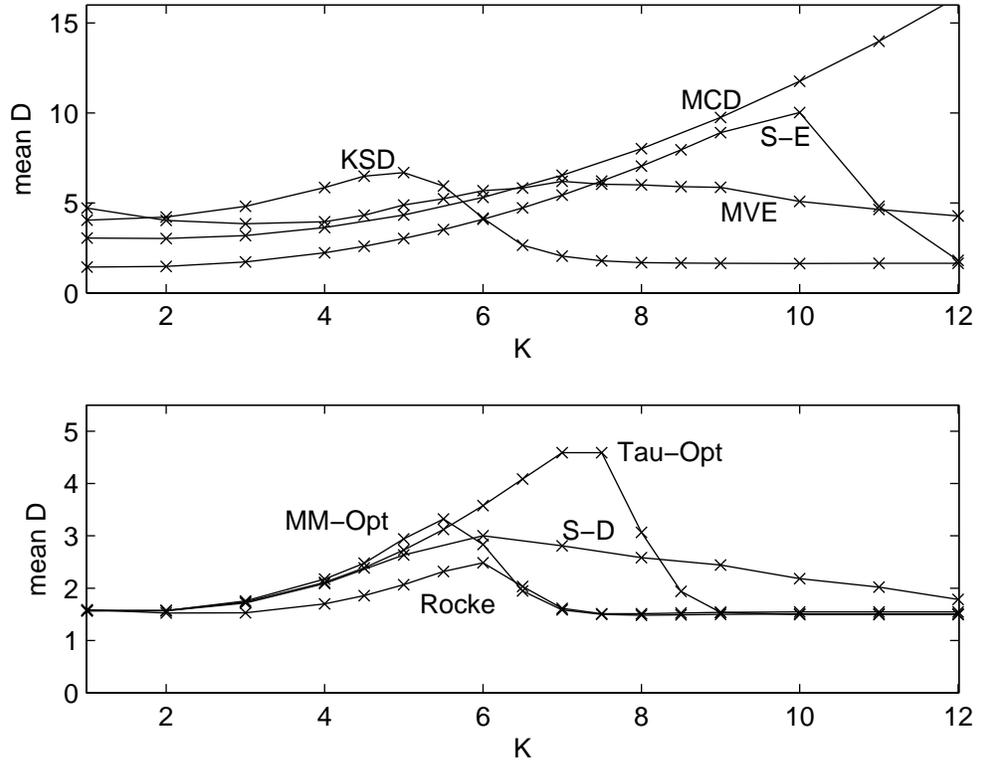}%
\caption{Mean $D$ of of scatter estimators for $p=20,$ $n=200,$ $\varepsilon
=0.1$ and $\gamma=0$ as a function of the outlier size $K.$ All iterative
estimators start from KSD.}%
\label{figConta10_20_0}%
\end{center}
\end{figure}

Figure \ref{figConta10_20_0} shows the values of $\overline{D}$ as a function
of the outlier size $K$ for some of the estimators in the case $p=20,$ $n=200$
and $\gamma=0.$ Here \textquotedblleft MM-Opt\textquotedblright\ stands for
\textquotedblleft MM with 'optimal' $\rho$\textquotedblright. All estimators
in the second panel start from KSD.%

The plot confirms the superiority of Rocke+KSD.

\subsection{Comparison with a non-equivariant estimator\label{SecDetermi}}

Recently Hubert et al, (2015) proposed two deterministic estimators, called
DetS and DetMM, of which the latter has a tuneable efficiency. We compare it
with Rocke+SD. The nominal efficiency of DetMM is chosen as 0.90. The scenario
is the same as above. However, since DetMM is not equivariant, the model is
now N$_{p}\left(  \mathbf{0,\Sigma}_{0}\right)  $ where $\mathbf{\Sigma}_{0}$
has unit diagonal elements and all non-diagonal elements equal to $\rho.$ We
chose the extreme cases $\rho=0$ and $=0.9.$ Since both yield qualitatively
similar results, we show in Table \ref{TabSimuDeter} only the results from the
first case.

\begin{center}%
\begin{table}[tbp] \centering
\begin{tabular}
[c]{cccccccccc}\hline
$\varepsilon$ & $\gamma$ &  & \multicolumn{3}{c}{Scatter} &  &
\multicolumn{3}{c}{Location}\\\hline
&  & $n=$ & 100 & 200 & 400 &  & 100 & 200 & 400\\\hline
0.1 & 0 & \multicolumn{1}{l}{Rocke+KSD} & \multicolumn{1}{r}{3.96} &
\multicolumn{1}{r}{2.41} & \multicolumn{1}{r}{1.63} & \multicolumn{1}{r}{} &
\multicolumn{1}{r}{0.39} & \multicolumn{1}{r}{0.27} & \multicolumn{1}{r}{0.19}%
\\
&  & \multicolumn{1}{l}{\ \ DetMM} & \multicolumn{1}{r}{26.97} &
\multicolumn{1}{r}{26.59} & \multicolumn{1}{r}{26.19} & \multicolumn{1}{r}{} &
\multicolumn{1}{r}{5.93} & \multicolumn{1}{r}{5.72} & \multicolumn{1}{r}{5.32}%
\\
& 0.5 & Rocke+KSD & \multicolumn{1}{r}{4.73} & \multicolumn{1}{r}{2.95} &
\multicolumn{1}{r}{2.38} & \multicolumn{1}{r}{} & \multicolumn{1}{r}{0.54} &
\multicolumn{1}{r}{0.41} & \multicolumn{1}{r}{0.36}\\
&  & DetMM & \multicolumn{1}{r}{18.01} & \multicolumn{1}{r}{17.89} &
\multicolumn{1}{r}{18.99} & \multicolumn{1}{r}{} & \multicolumn{1}{r}{3.46} &
\multicolumn{1}{r}{3.29} & \multicolumn{1}{r}{3.18}\\\hline
0.2 & 0 & Rocke+KSD & \multicolumn{1}{r}{10.62} & \multicolumn{1}{r}{5.22} &
\multicolumn{1}{r}{3.58} & \multicolumn{1}{r}{} & \multicolumn{1}{r}{1.47} &
\multicolumn{1}{r}{0.73} & \multicolumn{1}{r}{0.52}\\
&  & DetMM & \multicolumn{1}{r}{213.66} & \multicolumn{1}{r}{164.18} &
\multicolumn{1}{r}{156.82} & \multicolumn{1}{r}{} & \multicolumn{1}{r}{81.29}
& \multicolumn{1}{r}{78.19} & \multicolumn{1}{r}{77.50}\\
& 0.5 & Rocke+KSD & \multicolumn{1}{r}{12.08} & \multicolumn{1}{r}{9.24} &
\multicolumn{1}{r}{8.67} & \multicolumn{1}{r}{} & \multicolumn{1}{r}{2.33} &
\multicolumn{1}{r}{1.95} & \multicolumn{1}{r}{1.90}\\
&  & DetMM & \multicolumn{1}{r}{118.79} & \multicolumn{1}{r}{111.87} &
\multicolumn{1}{r}{109.79} & \multicolumn{1}{r}{} & \multicolumn{1}{r}{46.96}
& \multicolumn{1}{r}{46.90} & \multicolumn{1}{r}{45.81}\\\hline
\end{tabular}
\caption{Comaprison of Rocke and DetMM estimators: Maximum  mean D for $\rho=0$}\label{TabSimuDeter}%
\end{table}%

\end{center}

The performance of DetMM is clearly poor. We have not been able to find an
explanation for this disappointing behavior.

\subsection{Computing times\label{SecTimes}}

We compare the computing times of the Rocke estimator with MVE and KSD starts,
and of DetMM. The results are the average of 10 runs with normal samples, on a
PC with Intel TM12 Duo CPU and 3.01 GHz. The values of $n$ were $5p,$ $10p$
and $20p,$ with $p$ between 10 and 100. The number of subsamples for the MVE
was made to increase slowly as $50p.$ Table \ref{TabTiempos} displays the
results, where for brevity we show only the values for $p=20,$ 50, 80 and 100.

\begin{center}%
\begin{table}[tbp] \centering
\begin{tabular}
[c]{ccccc}\hline
$p$ & $n$ & Rocke+MVE & Rocke+KSD & DetMM\\\hline
\multicolumn{1}{r}{20} & \multicolumn{1}{r}{100} & \multicolumn{1}{r}{0.62} &
\multicolumn{1}{r}{0.06} & \multicolumn{1}{r}{0.20}\\
\multicolumn{1}{r}{} & \multicolumn{1}{r}{200} & \multicolumn{1}{r}{0.98} &
\multicolumn{1}{r}{0.079} & \multicolumn{1}{r}{0.30}\\
\multicolumn{1}{r}{} & \multicolumn{1}{r}{400} & \multicolumn{1}{r}{1.31} &
\multicolumn{1}{r}{0.15} & \multicolumn{1}{r}{0.54}\\
\multicolumn{1}{r}{50} & \multicolumn{1}{r}{250} & \multicolumn{1}{r}{5.03} &
\multicolumn{1}{r}{0.51} & \multicolumn{1}{r}{1.61}\\
\multicolumn{1}{r}{} & \multicolumn{1}{r}{500} & \multicolumn{1}{r}{6.54} &
\multicolumn{1}{r}{1.22} & \multicolumn{1}{r}{3.11}\\
\multicolumn{1}{r}{} & \multicolumn{1}{r}{1000} & \multicolumn{1}{r}{12.72} &
\multicolumn{1}{r}{3.07} & \multicolumn{1}{r}{6.51}\\
\multicolumn{1}{r}{80} & \multicolumn{1}{r}{400} & \multicolumn{1}{r}{14.55} &
\multicolumn{1}{r}{6.43} & \multicolumn{1}{r}{5.97}\\
\multicolumn{1}{r}{} & \multicolumn{1}{r}{800} & \multicolumn{1}{r}{22.46} &
\multicolumn{1}{r}{14.90} & \multicolumn{1}{r}{12.23}\\
\multicolumn{1}{r}{} & \multicolumn{1}{r}{1600} & \multicolumn{1}{r}{65.45} &
\multicolumn{1}{r}{22.48} & \multicolumn{1}{r}{26.64}\\
\multicolumn{1}{r}{100} & \multicolumn{1}{r}{500} & \multicolumn{1}{r}{26.86}
& \multicolumn{1}{r}{59.18} & \multicolumn{1}{r}{11.79}\\
\multicolumn{1}{r}{} & \multicolumn{1}{r}{1000} & \multicolumn{1}{r}{74.01} &
\multicolumn{1}{r}{91.63} & \multicolumn{1}{r}{24.47}\\
\multicolumn{1}{r}{} & \multicolumn{1}{r}{2000} & \multicolumn{1}{r}{152.06} &
\multicolumn{1}{r}{113.54} & \multicolumn{1}{r}{47.41}\\\hline
\end{tabular}
\caption{Mean computing times of estimators in seconds}\label{TabTiempos}%
\end{table}%

\end{center}

It is seen that Rocke+KSD is faster than DetMM for $p\leq80,$ and Rocke+MVE .
However it is slower than DetMM for $p=100.$ This rapid increase in computing
time is probably due to the optimization procedure employed by KSD, and may be
improved upon by choosing a more efficient optimizer.

\section{A real example\label{secReal}}

\begin{figure}
[tbh]
\begin{center}
\includegraphics[
height=4.4269in,
width=5.8878in
]%
{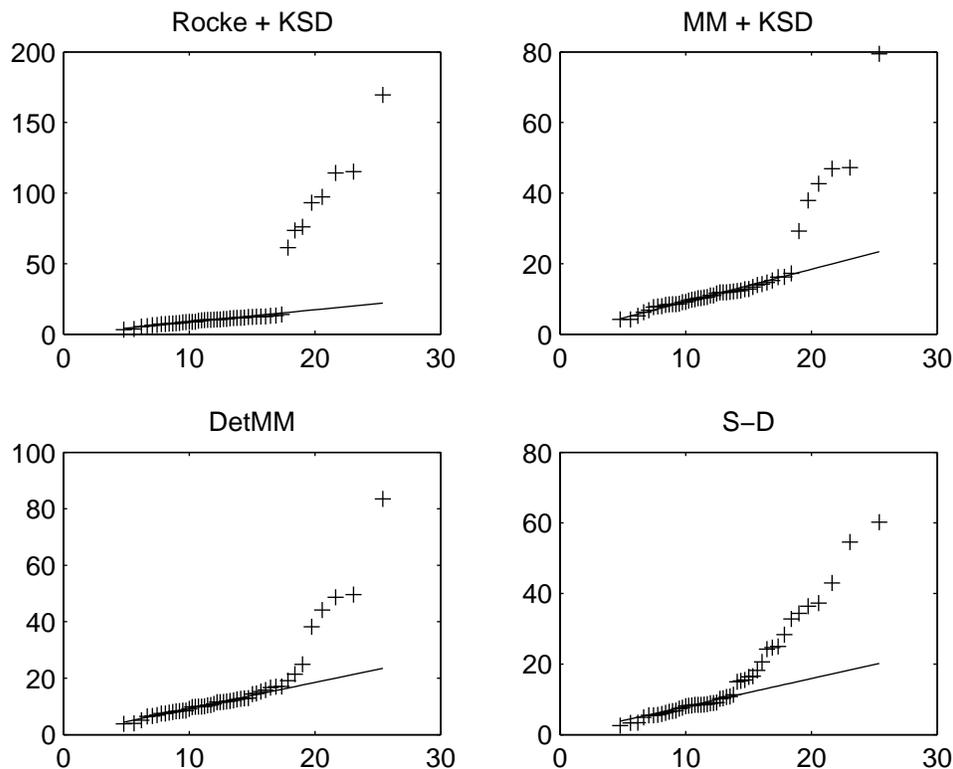}%
\caption{Wine data: Ordered Mahalanobis distances vs. $\chi_{p}^{2}%
$-quantiles.}%
\label{figWine1}%
\end{center}
\end{figure}

\begin{figure}
[tbh]
\begin{center}
\includegraphics[
height=4.4269in,
width=5.8878in
]%
{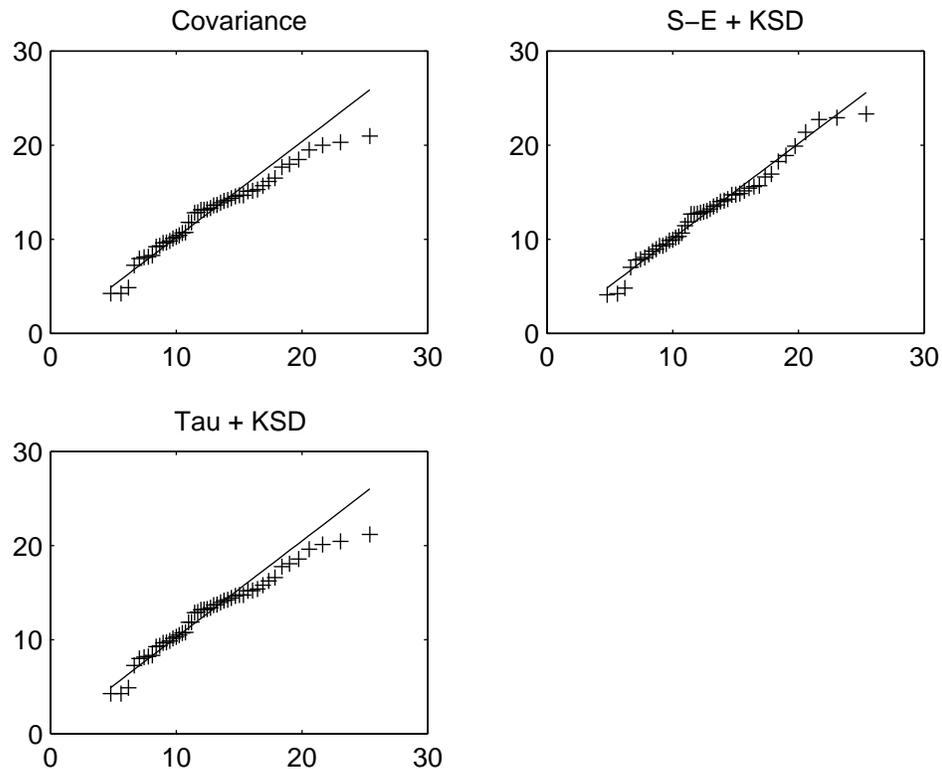}%
\caption{Wine data: Wine data: Ordered Mahalanobis distances vs. $\chi_{p}%
^{2}$-quantiles.}%
\label{figWine2}%
\end{center}
\end{figure}

We deal with the well-known wine data set, available at the UCI machine
learning repository: https://archive.ics.uci.edu/ml/datasets/Wine, which has
been employed as a benchmark data set for pattern recognition; see e. g.
(Aeberhard et al, 1994), and consists of three classes with 13 variables. The
estimators were applied to the data of class 3, with $n=48$ and $p=13.$ Since
KSD and MVE yielded similar results as initial estimators, we show only the
results corresponding to the former. Figures \ref{figWine1} and \ref{figWine2}
contain the QQ-plots of the (squared) Mahalanobis distances for the different estimators.%

%

Rocke, MM\ and DetMM pinpoint respectively 8, 6 and 5 possible outliers; S-D
seems to pinpoint an excessive number of possible outliers; while S-E and
$\tau$ behave like the classical estimator, showing no suspicious points. Some
subject-matter knowledge would be necessary to decide how atypical the
suspicious points are.

\section{Conclusions\label{secConclu}}

The Rocke estimator has a controllable efficiency for $p\geq15.$ With equal
efficiencies, the Rocke estimator with KSD start outperforms all its
competitors for shift contamination Its computing time is competitive for
$p<100,$ and can probably be improved upon. It can therefore be recommended
for estimation with $p\geq15$.

For $p<15$ we can recommend MM with \textquotedblleft
optimal\textquotedblright\ $\rho$ and KSD start.

\newpage

\section{Appendix\label{SecAppend}}

\subsection{Simulation results\label{secTablas}}

Tables \ref{TabSim05_10} to \ref{TabSim_3020} in this section contain the
detailed results of the simulation. In each scenario, the smallest and
next-to-smallest values are marked as bold and italic, respectively. This is
done only for estimators with controllable efficiency; in particular, Rocke is
not considered for $p\leq10,$ for its efficiency in this case is less than 0.90.

One would expect the values for a given estimator to decrease when $n$
increases. However, in many cases this does not hold for the estimators based
on the MVE. We have re-run the simulations with a different seed, and also
employed medians instead of means to rule out atypical cases, but this pattern
appears nevertheless. We have not been able to find an explanation for this
phenomenon. Since it always affects the largest values, it does not influence
the conclusions.

\begin{center}
\bigskip\bigskip%
\begin{table}[htb] \centering
\begin{tabular}
[c]{cccccccccc}\hline
$\gamma$ &  & Start & \multicolumn{3}{c}{Scatter} &  &
\multicolumn{3}{c}{Location}\\\hline
0 & $n=$ &  & 25 & 50 & 100 &  & 25 & 50 & 100\\\hline
\multicolumn{1}{r}{} & \multicolumn{1}{r}{S-E} & \multicolumn{1}{r}{MVE} &
\multicolumn{1}{r}{1.24} & \multicolumn{1}{r}{1.21} & \multicolumn{1}{r}{0.98}
& \multicolumn{1}{r}{} & \multicolumn{1}{r}{0.28} & \multicolumn{1}{r}{0.27} &
\multicolumn{1}{r}{0.24}\\
\multicolumn{1}{r}{} & \multicolumn{1}{r}{S-E} & \multicolumn{1}{r}{KSD} &
\multicolumn{1}{r}{1.27} & \multicolumn{1}{r}{1.09} & \multicolumn{1}{r}{0.90}
& \multicolumn{1}{r}{} & \multicolumn{1}{r}{0.28} & \multicolumn{1}{r}{0.27} &
\multicolumn{1}{r}{0.24}\\
\multicolumn{1}{r}{} & \multicolumn{1}{r}{MM-Bisq.} & \multicolumn{1}{r}{MVE}
& \multicolumn{1}{r}{1.26} & \multicolumn{1}{r}{1.28} &
\multicolumn{1}{r}{1.19} & \multicolumn{1}{r}{} &
\multicolumn{1}{r}{\textbf{0.28}} & \multicolumn{1}{r}{0.28} &
\multicolumn{1}{r}{0.24}\\
\multicolumn{1}{r}{} & \multicolumn{1}{r}{MM-Bisq.} & \multicolumn{1}{r}{KSD}
& \multicolumn{1}{r}{1.26} & \multicolumn{1}{r}{1.27} &
\multicolumn{1}{r}{1.24} & \multicolumn{1}{r}{} &
\multicolumn{1}{r}{\emph{0.29}} & \multicolumn{1}{r}{0.28} &
\multicolumn{1}{r}{0.25}\\
\multicolumn{1}{r}{} & \multicolumn{1}{r}{MM-Opt.} & \multicolumn{1}{r}{MVE} &
\multicolumn{1}{r}{1.29} & \multicolumn{1}{r}{0.87} & \multicolumn{1}{r}{0.74}
& \multicolumn{1}{r}{} & \multicolumn{1}{r}{0.29} & \multicolumn{1}{r}{0.24} &
\multicolumn{1}{r}{0.21}\\
\multicolumn{1}{r}{} & \multicolumn{1}{r}{MM-Opt.} & \multicolumn{1}{r}{KSD} &
\multicolumn{1}{r}{1.26} & \multicolumn{1}{r}{\textbf{0.85}} &
\multicolumn{1}{r}{\emph{0.72}} & \multicolumn{1}{r}{} &
\multicolumn{1}{r}{\emph{0.30}} & \multicolumn{1}{r}{0.24} &
\multicolumn{1}{r}{0.21}\\
\hline
\end{tabular}
\caption{Simulation:  Maximum mean $D$ for scatter and location with  $p=5$ and 10\%
contamination}\label{TabSim05_10}%
\end{table}%
%

\begin{table}[htb] \centering
\begin{tabular}
[c]{rccccccccc}\cline{2-10}%
$\gamma$ &  & Start & \multicolumn{3}{c}{Scatter} &  &
\multicolumn{3}{c}{Location}\\\cline{2-10}%
0 & $n=$ &  & 25 & 50 & 100 &  & 25 & 50 & 100\\\cline{2-10}
& \multicolumn{1}{r}{S-E} & \multicolumn{1}{r}{MVE} & \multicolumn{1}{r}{9.10}
& \multicolumn{1}{r}{6.02} & \multicolumn{1}{r}{6.23} & \multicolumn{1}{r}{} &
\multicolumn{1}{r}{2.63} & \multicolumn{1}{r}{2.15} & \multicolumn{1}{r}{2.34}%
\\
& \multicolumn{1}{r}{S-E} & \multicolumn{1}{r}{KSD} & \multicolumn{1}{r}{5.99}
& \multicolumn{1}{r}{4.38} & \multicolumn{1}{r}{4.26} & \multicolumn{1}{r}{} &
\multicolumn{1}{r}{1.93} & \multicolumn{1}{r}{1.70} & \multicolumn{1}{r}{1.87}%
\\
& \multicolumn{1}{r}{MM-Bisq.} & \multicolumn{1}{r}{MVE} &
\multicolumn{1}{r}{7.83} & \multicolumn{1}{r}{5.21} & \multicolumn{1}{r}{5.43}
& \multicolumn{1}{r}{} & \multicolumn{1}{r}{2.29} & \multicolumn{1}{r}{1.86} &
\multicolumn{1}{r}{2.03}\\
& \multicolumn{1}{r}{MM-Bisq.} & \multicolumn{1}{r}{KSD} &
\multicolumn{1}{r}{5.11} & \multicolumn{1}{r}{4.03} & \multicolumn{1}{r}{4.25}
& \multicolumn{1}{r}{} & \multicolumn{1}{r}{1.65} & \multicolumn{1}{r}{1.43} &
\multicolumn{1}{r}{1.57}\\
& \multicolumn{1}{r}{MM-Opti} & \multicolumn{1}{r}{MVE} &
\multicolumn{1}{r}{7.90} & \multicolumn{1}{r}{3.89} & \multicolumn{1}{r}{3.79}
& \multicolumn{1}{r}{} & \multicolumn{1}{r}{2.05} & \multicolumn{1}{r}{1.27} &
\multicolumn{1}{r}{1.31}\\
& \multicolumn{1}{r}{MM-Opti} & \multicolumn{1}{r}{KSD} &
\multicolumn{1}{r}{\textbf{3.90}} & \multicolumn{1}{r}{\textbf{2.27}} &
\multicolumn{1}{r}{\emph{2.05}} & \multicolumn{1}{r}{} &
\multicolumn{1}{r}{\textbf{1.17}} & \multicolumn{1}{r}{\textbf{0.83}} &
\multicolumn{1}{r}{\emph{0.80}}\\
\hline
\end{tabular}
\caption{Simulation:  Maximum mean $D$ for scatter and location with  $p=5$ and 20\%
contamination}\label{TabSim05_20}%
\end{table}%
\bigskip%

\begin{table}[htb] \centering
\begin{tabular}
[c]{cccccccccc}\hline
$\gamma$ &  & Start & \multicolumn{3}{c}{Scatter} &  &
\multicolumn{3}{c}{Location}\\\hline
0 & $n=$ &  & 50 & 100 & 200 &  & 50 & 100 & 200\\\hline
& S-E & MVE & \multicolumn{1}{r}{4.32} & \multicolumn{1}{r}{3.82} &
\multicolumn{1}{r}{3.59} & \multicolumn{1}{r}{} & \multicolumn{1}{r}{0.76} &
\multicolumn{1}{r}{0.69} & \multicolumn{1}{r}{0.67}\\
& S-E & KSD & \multicolumn{1}{r}{4.53} & \multicolumn{1}{r}{3.54} &
\multicolumn{1}{r}{3.21} & \multicolumn{1}{r}{} & \multicolumn{1}{r}{0.75} &
\multicolumn{1}{r}{0.64} & \multicolumn{1}{r}{0.62}\\
& MM-Bisq. & MVE & \multicolumn{1}{r}{3.56} & \multicolumn{1}{r}{3.23} &
\multicolumn{1}{r}{3.16} & \multicolumn{1}{r}{} & \multicolumn{1}{r}{0.66} &
\multicolumn{1}{r}{0.63} & \multicolumn{1}{r}{0.65}\\
& MM-Bisq. & KSD & \multicolumn{1}{r}{3.59} & \multicolumn{1}{r}{3.18} &
\multicolumn{1}{r}{2.98} & \multicolumn{1}{r}{} & \multicolumn{1}{r}{0.65} &
\multicolumn{1}{r}{0.62} & \multicolumn{1}{r}{0.60}\\
& MM-Opt. & MVE & \multicolumn{1}{r}{\emph{2.62}} & \multicolumn{1}{r}{1.78} &
\multicolumn{1}{r}{1.54} & \multicolumn{1}{r}{} &
\multicolumn{1}{r}{\emph{0.44}} & \multicolumn{1}{r}{0.34} &
\multicolumn{1}{r}{0.32}\\
& MM-Opt. & KSD & \multicolumn{1}{r}{\textbf{2.52}} &
\multicolumn{1}{r}{\emph{1.67}} & \multicolumn{1}{r}{\emph{1.45}} &
\multicolumn{1}{r}{} & \multicolumn{1}{r}{\textbf{0.43}} &
\multicolumn{1}{r}{\emph{0.33}} & \multicolumn{1}{r}{\emph{0.31}}\\
\hline
\end{tabular}
\caption{Simulation: Maximum mean $D$ for scatter and location with  $p=10$ and 10\%
contamination}\label{TabSim_1010}%
\end{table}%
%

\begin{table}[htb] \centering
\begin{tabular}
[c]{cccccccccc}\hline
$\gamma$ &  & Start & \multicolumn{3}{c}{Scatter} &  &
\multicolumn{3}{c}{Location}\\\hline
0 & $n=$ &  & 50 & 100 & 200 &  & 50 & 100 & 200\\\hline
& S-E & MVE & \multicolumn{1}{r}{15.20} & \multicolumn{1}{r}{12.83} &
\multicolumn{1}{r}{12.58} & \multicolumn{1}{r}{} & \multicolumn{1}{r}{5.12} &
\multicolumn{1}{r}{4.87} & \multicolumn{1}{r}{5.13}\\
& S-E & KSD & \multicolumn{1}{r}{13.65} & \multicolumn{1}{r}{11.26} &
\multicolumn{1}{r}{10.81} & \multicolumn{1}{r}{} & \multicolumn{1}{r}{4.2} &
\multicolumn{1}{r}{4.03} & \multicolumn{1}{r}{4.25}\\
& MM-Bisq. & MVE & \multicolumn{1}{r}{14.88} & \multicolumn{1}{r}{12.15} &
\multicolumn{1}{r}{11.86} & \multicolumn{1}{r}{} & \multicolumn{1}{r}{4.86} &
\multicolumn{1}{r}{4.58} & \multicolumn{1}{r}{4.76}\\
& MM-Bisq. & KSD & \multicolumn{1}{r}{9.69} & \multicolumn{1}{r}{8.01} &
\multicolumn{1}{r}{8.29} & \multicolumn{1}{r}{} & \multicolumn{1}{r}{3.31} &
\multicolumn{1}{r}{3.05} & \multicolumn{1}{r}{3.30}\\
& MM-Opt. & MVE & \multicolumn{1}{r}{11.67} & \multicolumn{1}{r}{8.72} &
\multicolumn{1}{r}{8.28} & \multicolumn{1}{r}{} & \multicolumn{1}{r}{3.01} &
\multicolumn{1}{r}{2.51} & \multicolumn{1}{r}{2.45}\\
& MM-Opt. & KSD & \multicolumn{1}{r}{\textbf{6.59}} &
\multicolumn{1}{r}{\textbf{3.88}} & \multicolumn{1}{r}{\textbf{3.52}} &
\multicolumn{1}{r}{} & \multicolumn{1}{r}{\textbf{1.83}} &
\multicolumn{1}{r}{\textbf{1.21}} & \multicolumn{1}{r}{\textbf{1.18}}\\
\hline
\end{tabular}
\caption{Simulation: Maximum mean $D$ for scatter and location with  $p=10$ and 20\%
contamination}\label{TabSim_10_20}%
\end{table}%
%

\begin{table}[htb] \centering
\begin{tabular}
[c]{cccccccccc}\hline
$\gamma$ &  & Start & \multicolumn{3}{c}{Scatter} &  &
\multicolumn{3}{c}{Location}\\\hline
& $n=$ &  & 75 & 150 & 300 &  & 75 & 150 & 300\\\hline
0 & S-E & MVE & 6.93 & 6.68 & 6.35 &  & 1.03 & 1.17 & 1.22\\
& S-E & KSD & 7.15 & 6.68 & 6.43 &  & 0.99 & 1.16 & 1.21\\
& MM-Bisq. & MVE & 5.47 & 5.65 & 5.65 &  & 0.88 & 1.13 & 1.23\\
& MM-Bisq. & KSD & 5.19 & 5.17 & 5.15 &  & 0.81 & 1.04 & 1.08\\
& MM-Opt. & MVE & 3.83 & 2.92 & 2.55 &  & 0.50 & 0.48 & 0.46\\
& MM-Opt. & KSD & \textbf{\ 3.33} & 2.38 & 2.25 &  & 0.46 & 0.38 & 0.41\\
\hline
\end{tabular}
\caption{Simulation: Maximum mean $D$ for scatter and location with  $p=15$ and 10\%
contamination}\label{TabSim_1510}%
\end{table}%
%

\begin{table}[htb] \centering
\begin{tabular}
[c]{cccccccccc}\hline
$\gamma$ &  & Start & \multicolumn{3}{c}{Scatter} &  &
\multicolumn{3}{c}{Location}\\\hline
& $n=$ &  & 75 & 150 & 300 &  & 75 & 150 & 300\\\hline
0 & S-E & MVE & \multicolumn{1}{r}{20.06} & \multicolumn{1}{r}{20.23} &
\multicolumn{1}{r}{20.43} & \multicolumn{1}{r}{} & \multicolumn{1}{r}{6.97} &
\multicolumn{1}{r}{7.91} & \multicolumn{1}{r}{8.44}\\
& S-E & KSD & \multicolumn{1}{r}{23.22} & \multicolumn{1}{r}{19.82} &
\multicolumn{1}{r}{19.59} & \multicolumn{1}{r}{} & \multicolumn{1}{r}{6.71} &
\multicolumn{1}{r}{6.77} & \multicolumn{1}{r}{7.3}\\
& MM-Bisq. & MVE & \multicolumn{1}{r}{14.06} & \multicolumn{1}{r}{21.25} &
\multicolumn{1}{r}{19.92} & \multicolumn{1}{r}{} & \multicolumn{1}{r}{8.95} &
\multicolumn{1}{r}{8.69} & \multicolumn{1}{r}{8.87}\\
& MM-Bisq. & KSD & \multicolumn{1}{r}{15.42} & \multicolumn{1}{r}{12.66} &
\multicolumn{1}{r}{12.15} & \multicolumn{1}{r}{} & \multicolumn{1}{r}{5.56} &
\multicolumn{1}{r}{4.99} & \multicolumn{1}{r}{5.02}\\
& MM-Opt & MVE & \multicolumn{1}{r}{10.75} & \multicolumn{1}{r}{13.87} &
\multicolumn{1}{r}{12.49} & \multicolumn{1}{r}{} & \multicolumn{1}{r}{2.54} &
\multicolumn{1}{r}{3.85} & \multicolumn{1}{r}{3.73}\\
& MM-Opt. & KSD & \multicolumn{1}{r}{10.41} & \multicolumn{1}{r}{\emph{5.68}}
& \multicolumn{1}{r}{\emph{5.10}} & \multicolumn{1}{r}{} &
\multicolumn{1}{r}{2.67} & \multicolumn{1}{r}{\emph{1.67}} &
\multicolumn{1}{r}{1.68}\\
\hline
\end{tabular}
\caption{Simulation: Maximum mean $D$ for scatter and location with  $p=15$ and 20\%
contamination}\label{TabSim_1520}%
\end{table}%
%

\begin{table}[htb] \centering
\begin{tabular}
[c]{cccccccccc}\hline
$\gamma$ &  & Start & \multicolumn{3}{c}{Scatter} &  &
\multicolumn{3}{c}{Location}\\\hline
& $n=$ &  & 100 & 200 & 400 &  & 100 & 200 & 400\\\hline
0 & S-E & MVE & \multicolumn{1}{r}{12.30} & \multicolumn{1}{r}{11.03} &
\multicolumn{1}{r}{10.69} & \multicolumn{1}{r}{} & \multicolumn{1}{r}{1.91} &
\multicolumn{1}{r}{1,79} & \multicolumn{1}{r}{1.71}\\
& S-E & KSD & \multicolumn{1}{r}{10.78} & \multicolumn{1}{r}{10.03} &
\multicolumn{1}{r}{10.06} & \multicolumn{1}{r}{} & \multicolumn{1}{r}{1.86} &
\multicolumn{1}{r}{1.72} & \multicolumn{1}{r}{1.64}\\
& MM-Bisq. & MVE & \multicolumn{1}{r}{8.82} & \multicolumn{1}{r}{8.14} &
\multicolumn{1}{r}{7.93} & \multicolumn{1}{r}{} & \multicolumn{1}{r}{1.83} &
\multicolumn{1}{r}{1.72} & \multicolumn{1}{r}{1.65}\\
& MM-Bisq. & KSD & \multicolumn{1}{r}{7.71} & \multicolumn{1}{r}{7.41} &
\multicolumn{1}{r}{7.36} & \multicolumn{1}{r}{} & \multicolumn{1}{r}{1.60} &
\multicolumn{1}{r}{1.59} & \multicolumn{1}{r}{1.50}\\
& MM-Opt. & MVE & \multicolumn{1}{r}{5.78} & \multicolumn{1}{r}{4.16} &
\multicolumn{1}{r}{3.52} & \multicolumn{1}{r}{} & \multicolumn{1}{r}{0.66} &
\multicolumn{1}{r}{0.59} & \multicolumn{1}{r}{0.56}\\
& MM-Opt. & KSD & \multicolumn{1}{r}{\emph{4.45}} & \multicolumn{1}{r}{3.32} &
\multicolumn{1}{r}{2.77} & \multicolumn{1}{r}{} & \multicolumn{1}{r}{0.54} &
\multicolumn{1}{r}{0.49} & \multicolumn{1}{r}{0.47}\\
\hline
\end{tabular}
\caption{Simulation: Maximum mean $D$ for scatter and location with  $p=20$ and 10\%
contamination}\label{TabSim20_10}%
\end{table}%

\end{center}

\bigskip

\begin{center}%
\begin{table}[htb] \centering
\begin{tabular}
[c]{cccccccccc}\hline
$\gamma$ &  & Start & \multicolumn{3}{c}{Scatter} &  &
\multicolumn{3}{c}{Location}\\\hline
& $n=$ &  & 100 & 200 & 400 &  & 100 & 200 & 400\\\hline
0 & S-E & MVE & \multicolumn{1}{r}{36.13} & \multicolumn{1}{r}{35.50} &
\multicolumn{1}{r}{32.45} & \multicolumn{1}{r}{} & \multicolumn{1}{r}{13.83} &
\multicolumn{1}{r}{12.71} & \multicolumn{1}{r}{11.60}\\
& S-E & KSD & \multicolumn{1}{r}{27.25} & \multicolumn{1}{r}{25.41} &
\multicolumn{1}{r}{25.38} & \multicolumn{1}{r}{} & \multicolumn{1}{r}{12.24} &
\multicolumn{1}{r}{11.95} & \multicolumn{1}{r}{10.64}\\
& MM-Bisq. & MVE & \multicolumn{1}{r}{40.38} & \multicolumn{1}{r}{38.61} &
\multicolumn{1}{r}{34.9} & \multicolumn{1}{r}{} & \multicolumn{1}{r}{17.34} &
\multicolumn{1}{r}{16.76} & \multicolumn{1}{r}{16.20}\\
& MM-Bisq. & KSD & \multicolumn{1}{r}{22.49} & \multicolumn{1}{r}{16.86} &
\multicolumn{1}{r}{16.18} & \multicolumn{1}{r}{} & \multicolumn{1}{r}{8.59} &
\multicolumn{1}{r}{7.57} & \multicolumn{1}{r}{7.38}\\
& MM-Opt. & MVE & \multicolumn{1}{r}{25.90} & \multicolumn{1}{r}{23.3} &
\multicolumn{1}{r}{20.33} & \multicolumn{1}{r}{} & \multicolumn{1}{r}{6.97} &
\multicolumn{1}{r}{6.33} & \multicolumn{1}{r}{5.92}\\
& MM-Opt. & KSD & \multicolumn{1}{r}{14.50} & \multicolumn{1}{r}{7.90} &
\multicolumn{1}{r}{7.49} & \multicolumn{1}{r}{} & \multicolumn{1}{r}{3.57} &
\multicolumn{1}{r}{2.22} & \multicolumn{1}{r}{2.38}\\
\hline
\end{tabular}
\caption{Simulation: Maximum mean $D$ for scatter and location with  $p=20$ and 20\%
contamination}\label{TabSim20_20}%
\end{table}%
%

\begin{table}[htb] \centering
\begin{tabular}
[c]{cccccccccc}\hline
$\gamma$ &  & Start & \multicolumn{3}{c}{Scatter} &  &
\multicolumn{3}{c}{Location}\\\hline
& $n=$ &  & 150 & 300 & 600 &  & 150 & 300 & 600\\\hline
0 & S-E & MVE & \multicolumn{1}{r}{17.90} & \multicolumn{1}{r}{17.03} &
\multicolumn{1}{r}{16.76} & \multicolumn{1}{r}{} & \multicolumn{1}{r}{2.64} &
\multicolumn{1}{r}{2.91} & \multicolumn{1}{r}{3.04}\\
& S-E & KSD & \multicolumn{1}{r}{20.13} & \multicolumn{1}{r}{18.39} &
\multicolumn{1}{r}{17.76} & \multicolumn{1}{r}{} & \multicolumn{1}{r}{2.63} &
\multicolumn{1}{r}{2.87} & \multicolumn{1}{r}{3.03}\\
& MM-Bisq. & MVE & \multicolumn{1}{r}{14.05} & \multicolumn{1}{r}{13.81} &
\multicolumn{1}{r}{13.59} & \multicolumn{1}{r}{} & \multicolumn{1}{r}{2.59} &
\multicolumn{1}{r}{3.24} & \multicolumn{1}{r}{3.51}\\
& MM-Bisq. & KSD & \multicolumn{1}{r}{12.02} & \multicolumn{1}{r}{11.80} &
\multicolumn{1}{r}{11.35} & \multicolumn{1}{r}{} & \multicolumn{1}{r}{2.64} &
\multicolumn{1}{r}{2.56} & \multicolumn{1}{r}{2.38}\\
& MM-Opt. & MVE & \multicolumn{1}{r}{8.89} & \multicolumn{1}{r}{7.44} &
\multicolumn{1}{r}{5.79} & \multicolumn{1}{r}{} & \multicolumn{1}{r}{0.87} &
\multicolumn{1}{r}{0.93} & \multicolumn{1}{r}{0.84}\\
& MM-Opt. & KSD & \multicolumn{1}{r}{7.03} & \multicolumn{1}{r}{5.34} &
\multicolumn{1}{r}{4.25} & \multicolumn{1}{r}{} & \multicolumn{1}{r}{0.73} &
\multicolumn{1}{r}{0.71} & \multicolumn{1}{r}{0.65}\\
\hline
\end{tabular}
\caption{Simulation: Maximum mean $D$ for scatter and location with  $p=30$ and 10\%
contamination}\label{TabSim_3010}%
\end{table}%
%

\begin{table}[htb] \centering
\begin{tabular}
[c]{cccccccccc}\hline
$\gamma$ &  & Start & \multicolumn{3}{c}{Scatter} &  &
\multicolumn{3}{c}{Location}\\\hline
& $n=$ &  & 150 & 300 & 600 &  & 150 & 300 & 600\\\hline
0 & S-E & MVE & \multicolumn{1}{r}{62.17} & \multicolumn{1}{r}{73.37} &
\multicolumn{1}{r}{76.03} & \multicolumn{1}{r}{} & \multicolumn{1}{r}{21.85} &
\multicolumn{1}{r}{28.44} & \multicolumn{1}{r}{30.79}\\
& S-E & KSD & \multicolumn{1}{r}{54.71} & \multicolumn{1}{r}{49.14} &
\multicolumn{1}{r}{47.32} & \multicolumn{1}{r}{} & \multicolumn{1}{r}{14.65} &
\multicolumn{1}{r}{15.51} & \multicolumn{1}{r}{16.19}\\
& MM-Bisq. & MVE & \multicolumn{1}{r}{77.51} & \multicolumn{1}{r}{110.69} &
\multicolumn{1}{r}{121.81} & \multicolumn{1}{r}{} & \multicolumn{1}{r}{31.98}
& \multicolumn{1}{r}{55.85} & \multicolumn{1}{r}{67.37}\\
& MM-Bisq. & KSD & \multicolumn{1}{r}{41.01} & \multicolumn{1}{r}{28.82} &
\multicolumn{1}{r}{27.49} & \multicolumn{1}{r}{} & \multicolumn{1}{r}{17.35} &
\multicolumn{1}{r}{15.02} & \multicolumn{1}{r}{14.08}\\
& MM-Opt. & MVE & \multicolumn{1}{r}{46.42} & \multicolumn{1}{r}{59.67} &
\multicolumn{1}{r}{65.08} & \multicolumn{1}{r}{} & \multicolumn{1}{r}{11.8} &
\multicolumn{1}{r}{16.91} & \multicolumn{1}{r}{19.82}\\
& MM-Opt. & KSD & \multicolumn{1}{r}{21.08} & \multicolumn{1}{r}{14.21} &
\multicolumn{1}{r}{14.54} & \multicolumn{1}{r}{} & \multicolumn{1}{r}{5.04} &
\multicolumn{1}{r}{3.79} & \multicolumn{1}{r}{4.33}\\
\hline
\end{tabular}
\caption{Simulation: Maximum mean $D$ for scatter and location with  $p=30$ and 20\%
contamination}\label{TabSim_3020}%
\end{table}%

\end{center}

\subsection{Approximations for the tuning constants\label{secTunRocke}}

For the MM estimator with MVE start, the constant $c$ in (\ref{defMM_min}) is
approximated by%
\[
c=\left(  a+\frac{b}{p}+\frac{c}{p^{2}}\right)  \left(  d+e\frac{p}{n}\right)
\]
where the coefficients depend on the $\rho-$function as follows:

\begin{center}%
\begin{tabular}
[c]{cccccc}\hline
$\rho$ & $a$ & $b$ & $c$ & $d$ & $e$\\\hline
\multicolumn{1}{l}{Bisquare} & \multicolumn{1}{r}{0.540} &
\multicolumn{1}{r}{3.538} & \multicolumn{1}{r}{-7.505} &
\multicolumn{1}{r}{1.114} & \multicolumn{1}{r}{-0.968}\\
\multicolumn{1}{l}{Optimal} & \multicolumn{1}{r}{0.469} &
\multicolumn{1}{r}{3.158} & \multicolumn{1}{r}{-0.928} &
\multicolumn{1}{r}{1.167} & \multicolumn{1}{r}{-1.698}\\\hline
\end{tabular}

\end{center}

For the MM estimator with KSD start, $c$ is approximated by%
\[
c=a+\frac{b}{p}+c\frac{p}{n},
\]
with

\begin{center}%
\begin{tabular}
[c]{cccc}\hline
$\rho$ & $a$ & $b$ & $c$\\\hline
\multicolumn{1}{l}{Bisquare} & \multicolumn{1}{r}{0.716} &
\multicolumn{1}{r}{2.572} & \multicolumn{1}{r}{-0.786}\\
\multicolumn{1}{l}{Optimal} & \multicolumn{1}{r}{0.612} &
\multicolumn{1}{r}{4.504} & \multicolumn{1}{r}{-1.112}\\\hline
\end{tabular}

\end{center}

\bigskip For the Rocke estimator, the value of $\alpha$ in (\ref{gama-alfa})
which yields 90\% efficiency is approximated for $p\geq15$ by%

\[
\alpha=ap^{b}n^{c}
\]
with%
\[%
\begin{tabular}
[c]{cccc}\hline
Start & $a$ & $b$ & $c$\\\hline
MVE & 0.00436 & -0.5030 & 0.4214\\
KSD & 0.00216 & -1.0078 & 0.8156\\\hline
\end{tabular}
\ \ \
\]

\bigskip

For the $\tau$-estimator the constant $c$ in (\ref{rho1-2}) which yields 90\%
efficiency is approximated for both MVE and KSD starts by%

\begin{equation}
c=ap^{b} \label{cpoten}%
\end{equation}

with

\begin{center}%
\begin{tabular}
[c]{ccc}\hline
$\rho$ & $a$ & $b$\\\hline
\multicolumn{1}{l}{Bisquare} & \multicolumn{1}{l}{6.2984} &
\multicolumn{1}{l}{-0.8458}\\
\multicolumn{1}{l}{Optimal} & \multicolumn{1}{l}{2.9987} &
\multicolumn{1}{l}{-0.4647}\\\hline
\end{tabular}

\end{center}

Finally, for the Stahel-Donoho estimator the weight function is $W\left(
r\right)  =W_{\mathrm{opt}}(d/c)$ where the values of $c$ that yield 90\%
efficiency are given by
\[
c=a+\frac{b}{n}+c\frac{p}{n},
\]
with%
\[%
\begin{tabular}
[c]{lccc}\hline
& $a$ & $b$ & $c$\\\hline
\multicolumn{1}{c}{Subsampling} & \multicolumn{1}{r}{5.116} & 63.820 &
\multicolumn{1}{r}{2.213}\\
\multicolumn{1}{c}{KSD} & \multicolumn{1}{r}{6.564} & 0.211 &
\multicolumn{1}{r}{24.286}\\\hline
\end{tabular}
\ \ \ \
\]

Here, $9c$ is the cutoff point(the value at which $W$ vanishes).

\subsection{The breakdown point of the KSD estimator\label{SecBP-KSD}}

The KSD procedure is defined in the same way as the Stahel-Donoho estimator,
but with a different set of directions. The population version is as follows.
Let $\mathbf{x}$ be a random vector with distribution $F.$ Let $U\subset
R^{p}$ be a set of directions $\mathbf{u}$ with $\left\Vert \mathbf{u}%
\right\Vert =1.$ Let $\mu$ and $\sigma$ be univariate robust location and
scale estimators. The outlyingness of a point $\mathbf{z\in}R^{p}$ is defined
as
\[
\mathrm{OL}\left(  \mathbf{z}\right)  =\max_{\mathbf{u\in}U}\frac
{|\mathbf{u}^{\prime}\mathbf{z-}\mu_{F}\left(  \mathbf{u}^{\prime}%
\mathbf{x}\right)  |}{\sigma_{F}\left(  \mathbf{u}^{\prime}\mathbf{x}\right)
}.
\]
The location and scale estimators are defined as weighted means and covariance
matrix with weights $W\left(  \mathrm{OL}\left(  \mathbf{x}\right)  \right)  $
where $W(t)\geq0$ is a nonincreasing function for $t\geq0.$

For a sample, the estimator is defined as above with $F$ the empirical
distribution. In the (theoretical) Stahel-Donoho estimator, $U$ is the set of
all directions, and $W$ is a smooth function; in actual practice, a finite set
of directions obtained by subsampling is employed.

The KSD procedure employs two sets of directions: $U=U_{1}\cup U_{2}$. The
first one is deterministic, and consists of a set of $p$ orthogonal directions
maximizing the kurtosis of $\mathbf{u}^{\prime}\mathbf{x}$ and $p $ directions
minimizing it. The other is a set of random \textquotedblleft specific
directions\textquotedblright\ obtained through a stratified sampling. We shall
deal only with the first one. Besides, $W$ is of \textquotedblleft hard
rejection\textquotedblright\ type:\ $W(t)=\mathbf{1}\left(  t\leq\beta\right)
$ where $\beta$ depends on $p.$

Theoretical calculations with KSD seem extremely difficult, and for this
reason we will limit ourselves to a very simplified case. We consider only the
population case with point-mass contamination; furthermore we assume that the
uncontaminated data are elliptically distributed. It will be shown that if
$\mu$ and $\sigma$ have breakdown 0.5, so has the KSD estimator.

Let $F_{0}$ be an elliptical distribution with fourth moments and consider the
contaminated distribution $F=\left(  1-\varepsilon\right)  F_{0}%
+\varepsilon\delta_{\mathbf{x}_{0}}$ with $\varepsilon<0.5.$ Because of the
estimator's equivariance it may be assumed that $F_{0}$ is radial, with zero
means and identity covariance matrix, and that $\mathbf{x}_{0}=K\mathbf{b}%
_{1}$ where $\mathbf{b}_{j}$ are the elements of the canonical base and $K>0$.
Put $A=\mathrm{E}_{F_{0}}x_{1}^{4},$ where $x_{1}$ is the first coordinate of
$\mathbf{x.}$ The rotational symmetry implies that $A=\mathrm{E}_{F_{0}%
}\left(  \mathbf{u}^{\prime}\mathbf{x}\right)  ^{4}$ for all $\mathbf{u}%
=\left(  u_{1},...,u_{p}\right)  ^{\prime}$ with $\left\Vert \mathbf{u}%
\right\Vert =1.$

We will show that the direction of the contamination, i.e., $\mathbf{u=b}%
_{1},$ is always included in the set. It is straightforward to show that the
kurtosis of a projection $\mathbf{u}^{\prime}\mathbf{x}$ under $F$ is%

\begin{equation}
\mathrm{Kurt}_{F}\left(  \mathbf{u}^{\prime}\mathbf{x}\right)  =g\left(
s\right)  =:\frac{a+bsK^{2}+cs^{2}K^{4}}{\left(  1+\varepsilon sK^{2}\right)
}, \label{kurto1}%
\end{equation}
where $s=u_{1}^{2}$ and%
\[
a=\left(  1-\varepsilon\right)  A,\ b=6\left(  1-\varepsilon\right)
\varepsilon^{2},\ c=\varepsilon\left(  1-4\varepsilon+6\varepsilon
^{2}-3\varepsilon^{3}\right)  \
\]

It follows that $\mathrm{Kurt}_{F}\left(  \mathbf{u}^{\prime}\mathbf{x}%
\right)  $ depends on $\mathbf{u}$ only through $s=u_{1}^{2}\in\lbrack0,1].$ A
laborious but straightforward calculation shows that the derivative of
$\ g\left(  s\right)  $ has the form $g^{\prime}\left(  s\right)  =u\left(
s\right)  v\left(  s\right)  ,$ where $u\left(  s\right)  >0$ does not depend
on $K,$ and
\[
v\left(  s\right)  =\left(  1-\varepsilon\right)  \left(  3\varepsilon
-A\right)  +sK^{2}\left(  1-4\varepsilon+3^{2}\right)  .
\]

The location of the extrema depends only on the sign of $v.$ Although the
result holds in general, to simplify the analysis we consider only the case
$A>1.5$. and we assume $K^{2}>A.$ There are two cases. If $\varepsilon
\geq1/3,$ then $v\left(  s\right)  <0$ for $s\in\lbrack0,1],$ and therefore
$\mathbf{u=b}_{1}$ is a minimizer of $\mathrm{Kurt}_{F}\left(  \mathbf{u}%
^{\prime}\mathbf{x}\right)  .$ If $\varepsilon<1/3$ there are maxima at $s=1$
and $s=0,$ and therefore the set of maximizing directions contains
$\mathbf{b}_{1}$ and a set of orthogonal $\mathbf{u}$'s which are orthogonal
to $\mathbf{b}_{1}$

It follows that
\[
\mathrm{OL}\left(  \mathbf{x}_{0}\right)  \geq\frac{|K-\mu\left(
x_{1}\right)  |}{\sigma\left(  x_{1}\right)  }.
\]
Note that $\mu\left(  x_{1}\right)  $ and $\sigma\left(  x_{1}\right)  $
depend on $K,$ but since $\varepsilon<0.5$ they are bounded. Therefore for $K$
large enough, $\mathrm{OL}\left(  \mathbf{x}_{0}\right)  $ will be larger than
the cutoff value $\beta$ and will therefore have null weight. This finishes
the proof.

\bigskip

\medskip

\textbf{References}

Aeberhard, S., Coomans, D. and De Vel, O. Comparative analysis of statistical
pattern recognition methods in high dimensional settings. Pattern Recognition
27,\textbf{ }1065-1077 (1994).

Agostinelli, C. Leung, A., Yohai, V.J. and Zamar, R.H.. Robust estimation of
multivariate location and scatter in the presence of cellwise and casewise
contamination. Test. DOI 10.1007/s11749-015-0453-3 (2015).

Croux, C. and Haesbroeck, G. Influence function and efficiency of the Minimum
Covariance Determinant scatter matrix estimator. J. Mult. Anal., 71, 161-190
(1999). 

Davies, P.L. Asymptotic behavior of S-estimates of multivariate location
parameters and dispersion matrices. Ann. Statist. 15\textbf{,} 1269--1292 (1987).

Donoho, D. L. Breakdown Properties of Multivariate Location Estimators,\ Ph.
D. Qualifying paper, Harvard University. (1982).

Filzmoser, P, Maronna, R. and Werner, M. Outlier identification in high
dimensions, Comp. Stat. \& Data Anal. 52, 1694--1711 (2008).

Hampel, F. R., Ronchetti, E. M., Rousseeuw, P. J., and Stahel, W. A.. Robust
Statistics: The Approach Based on Influence Functions. Wiley, New York (1986).

Hubert, M, Rousseeuw, P, Vanpaemela, D and Verdonck, T. The DetS and DetMM
estimators for multivariate location and scatter. Comput. Stat. \& Data Anal.
81\textbf{, } 64--75 (2015).

Kent, J.T. and Tyler, D.E. Constrained M-estimation for multivariate location
and scatter. Ann. Statist., 24, 1346-1370 (1996).

Lopuha\"{a}, H.P. Multivariate $\tau$-estimators for location and scatter.
Canad. J. Statist. 19, 307--321 (1991).

Lopuha\"{a}, H.P. Highly efficient estimators of multivariate location with
high breakdown point. Ann Stat., 20 398-413 (1992).

Maronna, R.A. Robust M-Estimators of Multivariate Location and Scatter. Ann.
Statist., 4 51-67 (1976).

Maronna, R.A., Martin, R.D. and Yohai, V.J. Robust Statistics: Theory and
Methods. John Wiley and Sons, New York (2006).

Maronna, R. A. and Yohai, V. J. The behavior of the Stahel--Donoho robust
multivariate estimator. J. Amer. Statist. Ass., 90, 330--341 (1995).

Maronna, R.A. and Yohai, V.J. High finite-sample efficiency and robustness
based on distance-constrained maximum likelihood. Comput. Stat. \& Data
Anal.\emph{.} DOI 10.1016/j.csda.2014.10.015.

Maronna, R and Zamar, R. Robust estimation of location and dispersion for
high-dimensional data sets. Technometrics 44, 307-317 (2002).

Muler, N. and Yohai, V.J. Robust estimates for ARCH processes. J. Time Ser.
Anal. 23\textbf{,} 341--375 (2002).

Paindaveine, D. and Van Bever, G. Inference on the shape of elliptical
distributions based on the MCD. J. Multivar. Anal. 129, 125-144 (2014).

Pe\~{n}a, D. and Prieto, F.J. Combining Random and Specific Directions for
Robust Estimation of High-Dimensional Multivariate Data. J. Comput. \& Graph.
Statist. 16, 228-254 (2007).

Rocke, D. Robustness properties of S-estimators of multivariate location and
shape in high dimension. \emph{\ }Ann. Statist. 24, 1327-1345 (1996).

Rousseeuw P.J. Multivariate estimation with high breakdown point. In:
Grossmann W., Pflug G., Vincze I., and Wertz W. (Eds.), Mathematical
Statistics and Applications, Vol. B, 283--297. Reidel, Dordrecht (1985).

Stahel, W.\ A.\ Breakdown of covariance estimators,\ Research report 31,
Fachgruppe f\"{u}r Statistik, E.T.H. Z\"{u}rich (1981),

Tatsuoka, K.S. and Tyler, D.E. On the uniqueness of S-functionals and
M-functionals under nonelliptical distributions. Ann. Statist. 28, 1219--1243 (2000).

Tyler, D. E. Finite-sample breakdown points of projection-based multivariate
location and scatter statistics. Ann. Statist., 22, 1024--1044 (1994).

Verboven, S. and Hubert, M.  LIBRA: a MATLAB library for robust analysis.
Chemometr. Intell. Lab., 75, 127--136 (2005).

Yohai, V.J. High Breakdown-Point and High Efficiency Robust Estimates for
Regression. Ann. Statist., 15, 642--656 (1987).

Yohai, V.J. and Zamar, R. Optimal locally robust M-estimates of regression. J.
Statist. Plan. \& Inference., 57, 73-92 (1997).

Yohai, V.J. and Zamar, R. High Breakdown-Point Estimates of Regression by
Means of the Minimization of an Efficient Scale. J. Amer. Statist. Assoc. 83,
406--413 (1988).

Zuo, Y., Cui, H. and He, X. On the Stahel-Donoho estimator and depth-weighted
means of multivariate data. Ann. Statist, 32, 167--188 (2004).

\medskip

\bigskip\medskip\textbf{Acknowledgements: }This work has been partially
supported by grants PICT 2011-0397 from ANPCyT and 20020130100279BA from
Universidad de Buenos Aires at Buenos Aires, Argentina.

\end{document}